\documentclass[11pt,leqno]{article}

\usepackage{amsmath,amsfonts,amscd,amssymb,theorem}

\long\def\comment#1\endcomment{}

\comment
\endcomment

\makeatletter
\begingroup
\gdef\th@dotted{\normalfont\itshape
  \def\@begintheorem##1##2{%
        \item[\hskip\labelsep \theorem@headerfont ##1\ ##2.]}%
\def\@opargbegintheorem##1##2##3{%
   \item[\hskip\labelsep \theorem@headerfont ##1\ ##2\ (##3).]}}
\endgroup
\makeatother

\theoremstyle{dotted}

\newtheorem{theorem}{Theorem}[section]
\newtheorem{lemma}[theorem]{Lemma}

\newtheorem{prop}[theorem]{Proposition}
\newtheorem{corr}[theorem]{Corollary}

\makeatletter
\begingroup
\gdef\th@upshape{\normalfont
  \def\@begintheorem##1##2{%
        \item[\hskip\labelsep \theorem@headerfont ##1\ ##2.]}%
\def\@opargbegintheorem##1##2##3{%
   \item[\hskip\labelsep \theorem@headerfont ##1\ ##2\ (##3).]}}
\endgroup
\makeatother

\theoremstyle{upshape}

\newtheorem{defn}[theorem]{Definition}
\newtheorem{remark}[theorem]{Remark}
\newtheorem{exa}[theorem]{Example}

\makeatletter
\renewcommand{\subsection}{\@startsection{subsection}{2}{0pt}{-3ex
plus -1ex minus -0.2ex}{-2mm plus -0pt minus
-2pt}{\normalfont\bfseries}} 
\renewcommand{\subsubsection}{\@startsection{subsubsection}{3}{0pt}{-3ex
plus -1ex minus -0.2ex}{-2mm plus -0pt minus
-2pt}{\normalfont\bfseries}} 
\makeatother

\makeatletter
\@addtoreset{equation}{section}
\makeatother

\newcommand{\cntrct}                % contraction with a vector field
{\hspace{2pt}\raisebox{1pt}{\text{$\lrcorner$}}\hspace{2pt}}

\newcommand{\proof}[1][Proof.]{\smallskip\noindent{\em #1}}
\def\endproof{\hfill\ensuremath{\square}\par\medskip}

\def\eqref#1{\thetag{\ref{#1}}}

\let\latexref=\ref
\def\ref#1{{\normalfont{\latexref{#1}}}}

\newcommand{\wt}{\widetilde}
\newcommand{\wh}{\widehat}

\newcommand{\dg}{\dagger}

\newcommand{\idot}{{\:\raisebox{1pt}{\text{\circle*{1.5}}}}}
\newcommand{\hdot}{{\:\raisebox{3pt}{\text{\circle*{1.5}}}}}
\newcommand{\eps}{\varepsilon}
\renewcommand{\phi}{\varphi}

\newcommand{\Fun}{\operatorname{Fun}}

\newcommand{\Hom}{\operatorname{Hom}}
\newcommand{\Ext}{\operatorname{Ext}}

\newcommand{\Ker}{\operatorname{Ker}}
\renewcommand{\Im}{\operatorname{Im}}

\newcommand{\gr}{\operatorname{\sf gr}}
\newcommand{\cchar}{\operatorname{\sf char}}
\newcommand{\ppt}{{\sf pt}}
\newcommand{\id}{\operatorname{\sf id}}

\newcommand{\D}{\mathcal{D}}
\newcommand{\DF}{\mathcal{D}F}

\newcommand{\Z}{\mathbb{Z}}
\newcommand{\Q}{\mathbb{Q}}

\newcommand{\K}{\mathbb{K}}

\newcommand{\lotimes}{\overset{\sf\scriptscriptstyle L}{\otimes}}

\newcommand{\C}{\mathcal{C}}

\newcommand{\I}{{\sf I}}

\newcommand{\vH}{\check{H}}
\newcommand{\vC}{\check{C}}

\newcommand{\Aut}{\operatorname{\sf Aut}}

\newcommand{\Exp}{\operatorname{Exp}}

\newcommand{\Per}{\operatorname{Per}}
\newcommand{\cPer}{\operatorname{\overline{Per}}}
\newcommand{\per}{\operatorname{per}}

\newcommand{\bCP}{\overline{CP}}
\newcommand{\bHP}{\overline{HP}}

\def\dlim_#1{{\displaystyle\lim_{#1}}^\hdot}

\newcommand{\E}{\mathcal{E}}
\newcommand{\hush}{\natural}

\newcommand{\amod}{{\text{\rm -mod}}}

\newcommand{\flfl}{{\flat\flat}}

\newcommand{\m}{\mathfrak{m}}

\title{Spectral sequences for cyclic homology}

\author{D. Kaledin\thanks{Partially supported by RSF, grant
    14-50-00005, and the Dynasty Foundation award}}

\date{\em To Maxim Kontsevich, for his 50th birthday.}

\begin{document}

\maketitle

\tableofcontents

\section*{Introduction.}

In modern non-commutative algebraic geometry as formalized for
example by M. Kontsevich and Y. Soibelman \cite{ks}, algebraic
varieties are replaced by associative unital DG algebras $A_\idot$
considered up to a derived Morita equivalence.  The role of
differential forms is played by Hochschild homology classes of
$A_\idot$, and de Rham cohomology corresponds to periodic cyclic
homology $HP_\idot(A_\idot)$. Periodic cyclic homology is related to
Hochschild homology by a standard spectral sequence. If $A_\idot$ is
derived Morita-equivalent to a smooth algebraic variety $X$, then
the spectral sequence reduces to the classical Hodge-to-de Rham
spectral sequence of P. Deligne. Because of this, one also calls it
the {\em Hodge-to-de Rham} spectral sequence in the general case.

If the smooth variety $X$ is also proper, and the base ring is a
field of characteristic $0$, then the Hodge-to-de Rham spectral
sequence degenerates. This is Deligne's reformulation of the
classical Hodge theory, and there is also an alternative purely
algebraic proof due to Deligne and L. Illusie \cite{di}.

Motivated by this, Kontsevich and Soibelman conjectured in \cite{ks}
that if a general DG algebra $A_\idot$ over a field of
characteristic $0$ is homologically smooth and homologically proper,
then the Hodge-to-de Rham spectral sequence degenerates.

\medskip

This conjecture has been largely proved in \cite{ka0}, by an
adaptation of the method of Deligne and Illusie. However, the
argument in that paper suffers from two related drawbacks:
\begin{itemize}
\item One has to impose an additional assumption that the DG algebra
  $A_\idot$ is question is concentrated in non-negative homological
  degrees (that is, $A_i=0$ for $i > 0$).
\item To pass form algebras to DG algebras, one uses Dold-Kan
  equivalence and simplicial methods. This is what forces one to
  impose the assumption above, and this is what makes a large part
  of the argument very hard to understand.
\end{itemize}
In this paper, we revisit the subject, and we give a proof of
Kontsevich-Soibelman Degeneration Conjecture free from any
additional technical assumptions.

\medskip

Our method is still the same --- it is based on reduction to
positive characteristic and adapting the approach of \cite{di}. In
principle, one could remove the assumption \thetag{$\bullet$} by
using simplicial-cosimplicial objects as e.g. in
\cite{ka6}. However, this would make the argument even more
opaque. Instead, we opt for an earlier approach tried in the
unpublished preprint \cite{ka-1}. This is much closer to \cite{di},
in that it uses explicitly two different spectral sequences that
exist in positive characteristic --- the Hodge-to-de Rham spectral
sequence, on one hand, and the so-called {\em conjugate spectral
  sequence} on the other hand. The spectral sequences are completely
different, but they have the same first page and the same last
page. So, what one actually proves is that the conjugate spectral
sequence degenerates, under some assumptions; the Hodge-to-de Rham
sequence then degenerates for dimension reasons.

For associative algebras of finite homological dimension, it is
relatively easy to construct a non-commutative version of the
conjugate spectral sequence, and this has been essentially done in
\cite{ka-1}. However, for general DG algebras, the question is much
more delicate. It took a while to realize that the sequence in
question simply does not exist --- or rather, it does exist, but
does not converge to periodic cyclic homology. What it converges to
is a completely new additive invariant of DG algebras and DG
categories constructed in \cite{ka4} under the name of {\em
  co-periodic cyclic homology} and denoted $\bHP_\idot(A_\idot)$. In
retrospect, this state of affairs has been also suggested by
Kontsevich 10 years ago \cite{Ko1,Ko2}, but the major push for
actually developing the theory has been given by recent works of
A. Beilinson \cite{be} and B. Bhatt \cite{bh}. Whatever the origins
of the theory are, now we know that co-periodic cyclic homology does
exist, it has a conjugate spectral sequence converging to it, and
for a homologically smooth and homologically bounded DG algebra
$A_\idot$, there is a comparison theorem providing a canonical
isomorphism $HP_\idot(A_\idot)\cong\bHP_\idot(A_\idot)$.  Therefore
a Deligne-Illusie type of argument for degeneration should still be
possible. This is what the present paper provides.

\bigskip

A couple of words about the organization of the paper. Out of
necessity, a large part of it is a recapitulation of my earlier
papers. In particular, Section~\ref{coper.sec} and
Section~\ref{conj.sec} contain the relevant results of
\cite{ka4}. Section~\ref{coper.sec} is concerned with general
results about co-periodic cyclic homology, summarized in
Theorem~\ref{per.thm}, and Section~\ref{conj.sec} contains the
construction of the conjugate spectral
sequence. Section~\ref{tate.sec} starts with some general results on
Tate (co)homology of finite groups. Everything is completely
standard but we do not know any good references (there are some
intersections with \cite[Subsection 6.3]{ka5}). Then we apply the
results to define a relative version of Tate cohomology, and use it
to prove results about co-periodic cyclic homology. In particular,
this includes a degeneration criterion
(Proposition~\ref{P.dege.prop}). At this point, we can already prove
all out degeneration results; however, we make a detour and use the
opportunity to correct one fault of \cite{ka4} --- namely, we
construct the conjugate spectral sequence in characteristic $2$, the
case excluded in \cite{ka4} for reasons explained in
\cite[Subsection 5.5]{ka4}. This is the subject of
Section~\ref{2.sec}. The technology used is a combination of some
splitting results of \cite[Section 6]{ka5} and the notion of a trace
functor of \cite{ka6}. Finally, in Section~\ref{dege.sec}, we prove
our degeneration results. This includes Theorem~\ref{hdr.thm}
equivalent to the Degeneration Conjecture of \cite{ks}.

\subsection*{Acknowledgement.} For several reasons at once, it is a
great pleasure to dedicate this paper to Maxim Kontsevich. It was
him who posed the original problem. He then actively encouraged my
earlier work on it, and at the same time, was consistently unhappy
about the unnecessary technical complications of \cite{ka0}. In
retrospect, it was also him who showed a better way to approach the
subject, although at the time, I didn't quite understand what he
meant. Another mathematician to whom this paper owes a lot is Vadim
Vologodsky. In particular, he always insisted that the notion of a
conjugate spectral sequence omitted in \cite{ka0} should be an
integral part of the story. Moreover, a crucial technical idea that
appears in \cite[Section 6]{ka5} is also due to him. I am also very
grateful to David Kazhdan and Vladimir Hinich for their interest in
this work, and for encouraging me to try to settle the question once
and for all. Finally, in the years that passed since \cite{ka0}, I
had an opportunity to discuss the subject with many people; I want
to specifically mention very useful conversations with Alexandre
Beilinson, Bhargav Bhatt, Boris Feigin, and Boris Tsygan.

\section{Co-periodic cyclic homology.}\label{coper.sec}

In this section, we recall main facts about co-periodic cyclic
homology introduced in \cite{ka4}, together with some terminology
and notation.

\subsection{Mixed complexes.}

A {\em mixed complex} $\langle V_\idot,B \rangle$ in an abelian
category $\E$ is a complex $V_\idot$ in $\E$ equipped with a map of
complexes $B:V_\idot \to V_\idot[-1]$ such that $B^2=0$ (we will
drop $B$ from notaton when it is clear from the context). The {\em
  periodic expansion} $\Per(V_\idot)$ of a mixed complex $V_\idot$
is the complex
$$
\Per(V_\idot) = V_\idot((u)),
$$
with the differential $d+Bu$, where $d$ is the differential in the
complex $V_\idot$, $u$ is a formal generator of cohomological degree
$2$, and $V_\idot((u))$ is shorthand for ``formal Laurent power
series in $u$ with coefficients in $V_\idot$''. Analogously, {\em
  co-periodic} and {\em polynomial periodic} expansions
$\cPer(V_\idot)$, $\per(V_\idot)$ are given by
$$
\cPer(V_\idot) = V_\idot((u^{-1})), \qquad \per(V_\idot) =
V_\idot[u,u^{-1}],
$$
again with the differential $d+uB$. By definition, the space of
Taylor power series $V_\idot[[u]] \subset V_\idot((u))$ is a
subcomplex in the periodic expansion $\Per(V_\idot)$; the {\em
  expansion} $\Exp(V_\idot)$ is the quotient complex
$$
\Exp(V_\idot) = V_\idot[u^{-1}] = V_\idot((u))/uV_\idot[[u]].
$$
Multiplication by $u$ induces an invertible degree-$2$ endomorphism of the
complexes $\Per(V_\idot)$, $\cPer(V_\idot)$, $\per(V_\idot)$ and a
non-invertible periodicity map $u:\Exp(V_\idot) \to
\Exp(V_\idot)[2]$. We have
\begin{equation}\label{cp.lim}
\Per(V_\idot) \cong \lim_{\overset{u}{\gets}}\Exp(V_\idot).
\end{equation}
Since a Laurent polynomial in $u$ is a Laurent power series both in
$u$ and in $u^{-1}$, we have natural functorial maps
\begin{equation}\label{V.per}
\begin{CD}
\Per(V_\idot) @<<< \per(V_\idot) @>>> \cPer(V_\idot)
\end{CD}
\end{equation}
for any mixed complex $V_\idot$. If $V_\idot$ is concentrated in a
finite range of degrees, then both maps are isomorphisms, but in
general, they are not.

\begin{exa}\label{cycl.mx.exa}
Assume given a module $E$ over a ring $R$, and assume that a cyclic
group $\Z/n\Z$ of some order $n$ acts on $E$, with $\sigma:E \to E$
being the generator. Then the length-$2$ complex
\begin{equation}\label{sigma.mx}
\begin{CD}
E @>{\id - \sigma}>> E
\end{CD}
\end{equation}
has a natural structure of a mixed complex, with the map $B$ given
by
$$
B = \id + \sigma + \dots + \sigma^{n-1}:E_\sigma \to E^\sigma.
$$
The expansion of the complex \eqref{sigma.mx} is the standard
homology complex $C_\idot(\Z/n\Z,E)$, and the periodic expansion is
the standard Tate homology complex $\vC_\idot(\Z/n\Z,E)$.
\end{exa}

\subsection{Small categories.}

In a sense, the standard example of a mixed complex that appears in
nature combines the complexes of Example~\ref{cycl.mx.exa} for all
integers $n \geq 1$. To package the data, it is convenient to use
the technology of homology of small categories. For any small
category $I$ and ring $R$, we denote by $\Fun(I,R)$ the abelian
category of functors from $I$ to $R$-modules, with the derived
category $\D(I,R)$, and for any functor $\gamma:I \to I'$ between
small categories, we denote by $\gamma^*:\Fun(I',R) \to \Fun(I,R)$
the pullback functor. The functor $\gamma^*$ is exact, and it has a
left and a right adjoint that we denote by
$\gamma_!,\gamma_*:\Fun(I,R) \to \Fun(I',R)$. For any $E \in
\Fun(I,R)$, the homology $H_\idot(I,E)$ of the category $I$ with
coefficients in $E$ is by definition given by
$$
H_i(I,E) = L^i\tau_!E,
$$
where $\tau:I \to \ppt$ is the tautological projection to the point
category $\ppt$.

The specific small category that we need is A. Connes' cyclic
category $\Lambda$ of \cite{C}. We do not reproduce here the full
definition (see e.g. \cite{lo}), but we do recall that objects of
$\Lambda$ correspond to cellular decompositions of the circle $S^1$,
and morphisms correspond to homotopy classes of cellular maps of a
certain type. We call $0$-cells {\em vertices}, and we call
$1$-cells {\em edges}. For any $n \geq 1$, there is exactly one
decomposition with $n$ vertices and $n$ edges. The corresponding
object in $\Lambda$ is denoted $[n]$, and we denote by $V([n])$ the
set of its vertices. Any map $f:[n'] \to [n]$ in $\Lambda$ induces a
map $f:V([n']) \to V([n])$. For any $v \in V([n])$, the preimage
$f^{-1}(v) \subset V([n'])$ carries a natural total order. We have
$\Aut([n])=\Z/n\Z$, the cyclic group, so that for any $E \in
\Fun(\Lambda,R)$, $E([n])$ is naturally an
$R[\Z/n\Z]$-module. Moreover, we have a natural embedding
$j:\Delta^o \to \Lambda$, where as usual, $\Delta$ is the category
of finite non-empty totally ordered sets, and $\Delta^o$ is the
opposite category. To keep notation consistent with the embedding
$j$, we denote by $[n] \in \Delta$ the set with $n$ elements.

By definition, the category $\Fun(\Delta^o,R)$ is category of
simplicial $R$-modu\-les, and for any $E \in \Fun(\Delta^o,R)$, we
have the standard chain complex $CH_\idot(E)$ with terms $CH_i(E) =
E([i+1])$, $i \geq 0$, and the differential $d_i:CH_i(E) \to
CH_{i-1}(E)$ given by
\begin{equation}\label{d.del}
d_i = \sum_{0 \leq j \leq i}\delta_j^i,
\end{equation}
where $\delta_\idot^\hdot$ are the face maps. Moreover, we also have
another complex $CH'_\idot(E)$ with the same terms as $CH_\idot(E)$,
but with the differential given by \eqref{d.del} with the summation
extended over $j$ from $0$ to $i-1$ (that is, we drop the last
term). Then the complex $CH_\idot(E)$ computes the homology
$H_\idot(\Delta^o,E)$, while the complex $CH'_\idot(E)$ is acyclic
--- in fact, canonically contractible.

Now, for any object $E \in \Fun(\Lambda,R)$, we have the simplicial
object $j^*E \in \Fun(\Delta^o,R)$, and it is well-known that the
complexes \eqref{sigma.mx} for $E([n])$, $n \geq 1$ fit together
into a single bicomplex
\begin{equation}\label{ch.mx}
\begin{CD}
CH'_\idot(E) @>{\id - \sigma^\dg}>> CH_\idot(E),
\end{CD}
\end{equation}
where for any $[n] \in \Lambda$, $\sigma \in \Aut([n])$ is the
generator of the cyclic group $\Z/n\Z$, and $\sigma^\dg =
(-1)^{n+1}\sigma$. We denote by $CH_\idot(E)$ the total complex of
the bicomplex \eqref{ch.mx}. Furthermore, it is also well-known that
the maps $B$ of Example~\ref{cycl.mx.exa} fit together into a single
map $B:CH_\idot(E) \to CH_\idot(E)[-1]$ that turns $CH_\idot(E)$
into a mixed complex.

It will be convenient to recast this construction in a slightly
different way. For any $[n] \in \Lambda$, denote by $\K_\idot([n])$
the standard celluar chain complex computing the homology of the
circle $S^1$ with respect to the decomposition correcponding to
$[n]$. Then it turns out that $K_\idot$ is functorial with respect
to morphisms in $\Lambda$, so that we obtain an exact sequence
\begin{equation}\label{4term}
\begin{CD}
0 @>>> \Z @>{\kappa_1}>> \K_1 @>>> \K_0 @>>> \Z @>{\kappa_0}>> 0
\end{CD}
\end{equation}
in the category $\Fun(\Lambda,\Z)$. For any $E \in \Fun(\Lambda,R)$,
denote $\K_\idot(E) = \K_\idot \otimes E$. Then $\K_\idot(E)$ is a
mixed complex in $\Fun(\Lambda,R)$, with the map $B$ given by
\begin{equation}\label{B.kappa}
B = (\kappa_1 \circ \kappa_0) \otimes \id.
\end{equation}
Now for any $E \in \Fun(\Lambda,R)$, denote by $cc_\idot(E)$ the
cokernel of the natural map \eqref{ch.mx}. Then one can show that we
have natural identifications
$$
CH_\idot(j^*E) \cong cc_\idot(\K_0(E)), \qquad CH'_\idot(j^*E) \cong
cc_\idot(\K_1(E)),
$$
and the mixed complex $CH_\idot(E)$ is then given by
\begin{equation}\label{ch.cc}
CH_\idot(E) \cong cc_\idot(\K_\idot(E)),
\end{equation}
with the map $B$ induced by the map $B$ of \eqref{B.kappa}.

Finally, observe that if we are given a complex $E_\idot$ in the
category $\Fun(\Lambda,R)$, then we can apply $cc_\idot(-)$ and
$CH_\idot(-)$ to $E_\idot$ termwise. We denote by
$cc_\idot(E_\idot)$, $CH_\idot(E_\idot)$ the sum-total complexes of
the resulting bicomplexes. Explicitly, we have
\begin{equation}\label{cc.exp}
cc_j(E_\idot) = \bigoplus_{n \geq 1}E_{j+n}([n])_{\sigma^\dg},
\qquad j \in \Z,
\end{equation}
with the differential induced by the differential \eqref{d.del} and
the differential in the complex $E_\idot$.

\subsection{Cyclic homology.}\label{cycl.subs}

We can now define periodic and co-perio\-dic cyclic homology.

\begin{defn}
Assume given a ring $R$ and a complex $E_\idot$ in the category
$\Fun(\Lambda,R)$. Then the {\em cyclic homology complex}
$CC_\idot(E_\idot)$, the {\em periodic cyclic homology complex}
$CP_\idot(E_\idot)$, the {\em co-periodic cyclic homology complex}
$\bCP_\idot(E_\idot)$, and the {\em polynomial periodic cyclic
  homology complex} $cp_\idot(E_\idot)$ are given by
\begin{equation}\label{b.cp}
\begin{aligned}
CC_\idot(E_\idot) &= \Exp(CH_\idot(E_\idot)), \qquad\: 
CP_\idot(E_\idot) = \Per(CH_\idot(E_\idot)),\\
\bCP_\idot(E_\idot) &= \cPer(CH_\idot(E_\idot)), \qquad\quad
cp_\idot(E_\idot) = \per(CH_\idot(E_\idot)).
\end{aligned}
\end{equation}
The {\em periodic} resp.\ {\em co-periodic} cyclic homology
$HP_\idot(E_\idot)$ resp.\ $\bHP_\idot(E_\idot)$ is the homology of
the complexes $CP_\idot(E_\idot)$ resp.\ $\bCP_\idot(E_\idot)$.
\end{defn}

We note that the first line in \eqref{b.cp} is completely standard;
it is the second line that defines new theories introduced in
\cite{ka4}. The homology of the complex $CC_\idot(E_\idot)$ is the
usual cyclic homology $HC_\idot(E_\idot)$, and it is well-known that
we have a natural identification
\begin{equation}\label{hc.hla}
HC_\idot(E_\idot) \cong H_\idot(\Lambda,E_\idot),
\end{equation}
where $E_\idot$ in the right-hand side is understood as the
corresponding object in the derived category $\D(\Lambda,R)$. One
can combine this with \eqref{cp.lim} to express
$HP_\idot(-)$. Co-periodic cyclic homology functor $\bHP_\idot(-)$
does not admit such a homological expression, and in fact, it is not
true that quasiisomorphic complexes have isomorphic
$\bHP_\idot$. For any complex $E_\idot$, we do have functorial maps
\begin{equation}\label{cp.per}
\begin{CD}
CP_\idot(E_\idot) @<<< cp_\idot(E_\idot) @>>> \bCP_\idot(E_\idot)
\end{CD}
\end{equation}
induced by the maps \eqref{V.per}, but in general, these maps are
not quasiisomorphisms. We also note that we have a natural
functorial map
\begin{equation}\label{CC.cc}
\alpha:CC_\idot(E_\idot) \to cc_\idot(E_\idot),
\end{equation}
and in general, this map is not a quasiisomorphism either. One
example where it is a quasiisomorphism is $E_\idot =
\K_\idot(E'_\idot)$ for some complex $E'_\idot$ in
$\Fun(\Lambda,R)$. In this case, by \cite[Lemma 3.11]{ka4}, $\alpha$
induces a natural isomorphism
\begin{equation}\label{hc.hh}
HC_\idot(\K_\idot(E'_\idot)) \cong HH_\idot(E'_\idot),
\end{equation}
where $HH_\idot(E'_\idot)$ is the homology of the complex
$CH_\idot(E'_\idot)$ (or equivalently, of the complex
$CH_\idot(j^*E'_\idot)$).

Assume now given a Noetherian commutative ring $k$ and a DG algebra
$A_\idot$ termwise-flat over $k$. Then one defines a complex
$A^\hush_\idot$ in $\Fun(\Lambda,k)$ as follows. For any $[n] \in
\Lambda$, $A^\hush_\idot([n]) \cong A_\idot^{\otimes_k n}$, with
terms of the product numbered by vertices $v \in V([n])$. For any
map $f:[n'] \to [n]$, the corresponding map $A^\hush_\idot(f)$ is
given by
\begin{equation}\label{a.hush}
A^\hush_\idot(f) = \bigotimes_{v \in V([n])}m_{f^{-1}(v)},
\end{equation}
where for any finite set $S$, $m_S:A^{\otimes_k S}_\idot \to
A_\idot$ is the multiplication map in the DG algebra $A_\idot$. Then
by definition, cyclic homology and periodic cyclic homology of the
DG algebra $A_\idot$ are the cyclic homology and the periodic cyclic
homology of the complex $A^\hush_\idot$, and we define co-periodic
and polynomial periodic cyclic homology $\bHP_\idot(A_\idot)$,
$hp_\idot(A_\idot)$ in the same way: we set
$$
\bHP_\idot(A_\idot) = \bHP_\idot(A^\hush_\idot), \qquad
hp_\idot(A_\idot) = hp_\idot(A^\hush_\idot).
$$
Here are, then, the main two results about $\bHP_\idot(A_\idot)$
proved in \cite{ka4}.

\begin{theorem}\label{per.thm}
\begin{enumerate}
\item Co-periodic cyclic homology functor $\bHP_\idot(-)$ extends to
  an additive invariant of small DG categories (in particular, a
  quasiisomorphism fo DG algebras induces an isomorphism of their
  co-peri\-o\-dic cyclic homology groups).
\item Assume that $k \otimes \Q = 0$ and $A_\idot$ is homologically
  smooth and homologically bounded over $k$. Then the maps
  \eqref{cp.per} for $E_\idot=A_\idot^\hush$ are quasiisomorphisms,
  so that we have
$$
HP_\idot(A_\idot) \cong hp_\idot(A_\idot) \cong \bHP_\idot(A_\idot).
$$
\end{enumerate}
\end{theorem}

\proof{} The first statement \thetag{i} is \cite[Theorem 6.6]{ka4},
and the second statement \thetag{ii} is contained in \cite[Theorem
  6.7]{ka4}.
\endproof

In Theorem~\ref{per.thm}~\thetag{ii}, {\em homologically smooth} as
usual means that $A_\idot$ is perfect as an $A_\idot$-bimodule (that
is, $A_\idot^o \otimes_k A_\idot$-module). If $k$ is a field, then
{\em homologically bounded} simply means that the complex $A_\idot$
has a finite number of non-trivial homology groups. In the general
case, this must hold for $A_\idot \otimes \overline{k}$, where
$\overline{k}$ is any residue field of the ring $k$.

\section{Conjugate spectral sequence.}\label{conj.sec}

For any mixed complex $\langle V_\idot,B \rangle$ in an abelian
category $\E$, the $u$-adic filtration on $\Per(V_\idot) \cong
V_\idot((u))$ induces a convergent spectral sequence
$$
H_\idot(V_\idot)((u)) \Rightarrow H_\idot(\Per(V_\idot)),
$$
where $H_\idot(-)$ stands for homology objects. In particular, for any
ring $R$ and complex $E_\idot \in \Fun(\Lambda,R)$, we have a
convergent spectral sequence
\begin{equation}\label{hh.hp.la}
HH_\idot(E_\idot)((u)) \Rightarrow HP_\idot(E_\idot).
\end{equation}
If $E_\idot = A^\hush_\idot$ for a DG algebra $A_\idot$ over a
commutative ring $k$, $HH_\idot(A^\hush_\idot)$ is naturally
identified with the Hochschild homology $HH_\idot(A_\idot)$ of the
DG algebra $A_\idot$, so that \eqref{hh.hp.la} reads as
\begin{equation}\label{hdr.sp}
HH_\idot(A_\idot)((u)) \Rightarrow HP_\idot(A_\idot).
\end{equation}
This is the {\em Hodge-to-de Rham spectral sequence}.

For co-periodic cyclic homology, no analog of \eqref{hh.hp.la} is
currently known, but under some assumptions, we do have a version of
\eqref{hdr.sp}. This was inroduced in \cite{ka4} under the name of
the {\em conjugate spectral sequence}. In this section, we briefly
recall the construction.

\subsection{Filtrations.}

The main technical tool used in \cite{ka4} for studying co-periodic
cyclic homology is the use of filtrations and filtered derived
categories. Filtrations are decreasing and indexed by all integers
--- that is, a {\em filtered complex} in an abelian category $\E$ is
a complex $E_\idot$ equipped with a collection of subcomplexes
$F^iE_\idot \subset E_\idot$, $i \in \Z$ such that $F^iE_\idot
\subset F^jE_\idot$ for $i \geq j$. A filtration $F^\hdot$ is {\em
  termwise-split} if for any $i$ and $j$, the embedding $F^iE_j \to
E_j$ admits a one-sided inverse $E_j \to F^iE_j$. The stupid
filtration $F^\hdot$ on a complex $E_\idot$ is obtained by setting
$F^jE_i=E_i$ if $i + j \leq 0$ and $0$ otherwise; it is
tautologically termwise-split. For any filtration $F^\hdot$ and any
integer $n \geq 1$, the {\em $n$-th rescaling} $F^\hdot_{[n]}$ of
$F^\hdot$ is given by $F^i_{[n]} = F^{in}$, $i \in \Z$, and the {\em
  shift by $n$} $F^\hdot_n$ is given by $F^i_n = F^{i+n}$. A map
$\langle E'_\idot,F'_\idot \rangle \to \langle E_\idot,F_\idot
\rangle$ of filtered complexes is a {\em filtered quasiisomorphism}
if for any $i$, the induced map $\gr^i_{F'}E'_\idot \to
\gr^i_FE_\idot$ of associated graded quotients is a quasiisomorphism
(there is no requirement on the map $E'_\idot \to E_\idot$ nor on
the maps $F^iE'_\idot \to F^iE_\idot$). Inverting filtered
quasiisomorphisms in the category of filtered complexes, we obtain
the {\em filtered derived category} $\DF(\E)$.

The {\em completion} $\wh{E}_\idot$ of a filtered complex $\langle
E_\idot,F^\hdot \rangle$ is given by
$$
\wh{E}_\idot =
\lim_{\overset{i}{\gets}}\lim_{\overset{j}{\to}}F^jE_\idot/F^iE_\idot \cong
\lim_{\overset{j}{\to}}\lim_{\overset{i}{\gets}}F^jE_\idot/F^iE_\idot,
$$
where the limit is over all integers $i \geq j$, with $i$ going to
$\infty$ and $j$ to $-\infty$. If a map $E'_\idot \to E_\idot$ is a
filtered quasiisomorphism, the induced map $\wh{E}'_\idot \to
\wh{E}_\idot$ of completions is a quasiisomorphism. The converse is
not true: two very different filtrations can have the same
completion. Specifically, say that two filtrations $F_1^\hdot$,
$F^\hdot_2$ on the same complex $E_\idot$ are {\em commensurable} if
for any integers $i$, $j$, there exist integers $j_1 \leq j \leq
j_2$ such that $F^{j_2}_2E_i \subset F^j_1E_i \subset F^{j_1}_2E_i$
and $F^{j_2}_1E_i \subset F^j_2E_i \subset F^{j_1}_1E_i$. Then two
commensurable filtrations on a complex $E_\idot$ obviously give the
same completion $\wh{E}_\idot$.

Every shift $F^\hdot_{n}$ and every rescaling $F^\hdot_{[n]}$ of a
filtration $F^\hdot$ is obviously commensurable to $F^\hdot$. For a
more non-trivial example, assume given a filtered complex $\langle
E_\idot,F^\hdot \rangle$, and define its {\em filtered truncations}
by setting
\begin{equation}\label{tau.F}
\begin{aligned}
\tau^F_{\geq n}E_i &= d^{-1}(F^{n+1-i}E_{i-1}) \cap F^{n-i}E_i
\subset E_i,\\
\tau^F_{\leq n}E_i &= E_i/(F^{n+1-i}E_i + d(F^{n-i}E_{i+1}))
\end{aligned}
\end{equation}
for any integer $n$. Denote also $\tau^F_{[n,m]} = \tau^F_{\geq
  n}\tau^F_{\leq m}$ for any integers $n \leq m$. Then the
subcomplexes $\tau^F_{\geq n}E_\idot \subset E_\idot$ for varying
$n$ give a filtration $\tau^\hdot$ on $E_\idot$, and this filtration
is commensurable to $F^\hdot$.

\begin{exa}
Assume that $E_\idot$ is the sum-total complex of a bicomplex
$E_{\idot,\idot}$, and let $F^\hdot$ be the stupid filtration with
respect to the first coordinate. Then $\tau^F_{\geq n}E_\idot$,
$\tau^F_{\leq n}E_\idot$ are canonical truncations with respect to
the second coordinate.
\end{exa}

In the general case, the truncation functors $\tau^F$ are also
related to the canonical filtrations: for any integers $i$, $j$, we
have natural isomorphisms
$$
\gr^j_F\tau^F_{\geq i}E_\idot \cong \tau_{\geq i-j}\gr^j_FE_\idot,
\qquad \gr^j_F\tau^F_{\leq i}E_\idot \cong \tau_{\leq
  i-j}\gr^j_FE_\idot,
$$
where $\tau_{\geq n}$, $\tau_{\leq n}$ are the usual canonical
truncations. In particular, the functors $\tau^F$ preserve filtered
quasiisomorphisms and descend to the filtered derived category
$\DF(\E)$, where they become truncation functors with respect to the
well-known $t$-structure of \cite{BBD} (the heart of this
$t$-structure is the category of complexes in $\E$). However, it
will be useful to have the truncation functors already on the level
of filtered complexes.

We note that since the filtrations $F^\hdot$ and $\tau^\hdot$ are
commensurable, $\gr^i_\tau E_\idot$ is complete with respect to
$F^\hdot$ for any integer $i$, so that for any filtered
quaiisomorphism $f:\langle E'_\idot,F'_\idot \rangle \to \langle
E_\idot,F_\idot \rangle$ is also a filtered quasiisomorphism with
respect to the filtrations $\tau$. Thus sending $\langle
E_\idot,F^\hdot \rangle$ to $\langle E_\idot,\tau^\hdot \rangle$
descends to an endofunctor
\begin{equation}\label{F.tau}
\DF(\E) \to \DF(\E)
\end{equation}
of the filtered derived category $\DF(\E)$.  We also note that for
any integer $i$ and filtered complex $\langle E_\idot,F^\hdot
\rangle$, $\gr^i_\tau E_\idot$ is quasiisomorphic (but not
isomorphic) to $\tau_{[i,i]}^FE_\idot$.

\medskip

Now, for any ring $R$ and any complex $E_\idot$ in the category
$\Fun(\Lambda,R)$ equipped with a termwise-split filtration
$F^\hdot$, define the {\em standard filtration} on the complex
$cc_\idot(E_\idot)$ by setting
\begin{equation}\label{st.cc}
F^icc_j(E_\idot) = \bigoplus_{n \geq
  1}F^{i+n}E_{j+n}([n])_{\sigma^\dg}, \qquad i,j \in \Z,
\end{equation}
where we use the decomposition \eqref{cc.exp} of the complex
$cc_\idot(E_\idot)$. By virtue of \eqref{ch.cc}, the standard
filtration extends to the complex $CH_\idot(E_\idot)$ and then to
its periodic expansions $\bCP_\idot(E_\idot)$,
$cp_\idot(E_\idot)$.

\begin{lemma}[{{\cite[Lemma 3.8]{ka4}}}]\label{cp.bcp.le}
Equip a complex $E_\idot$ in $\Fun(\Lambda,R)$ with the stupid
filtration $F^\hdot$. Then the co-periodic complex
$\bCP_\idot(E_\idot)$ is isomorphic to the completion of the
polynomial periodic complex $cp_\idot(E_\idot)$ with respect to the
standard filtration.\endproof
\end{lemma}

By virtue of this result, one can reduce the study of
$\bCP_\idot(E_\idot)$ to the study of $cp_\idot(E_\idot)$ equipped
with the standard filtration.

Let us now do the following. Fix an integer $p \geq 1$, and for any
complex $E_\idot$ in $\Fun(\Lambda,R)$ equipped with the stupid
filtration, denote by $cp_\idot(E_\idot)^{[p]}$ the complex
$cp_\idot(E_\idot)$ equipped with the $p$-th rescaling of the
standard filtration of \eqref{st.cc}.

\begin{defn}
The {\em conjugate filtration} $V^\hdot$ on $cp_\idot(E_\idot)$ is
given by
\begin{equation}\label{conj.glo}
V^ncp_\idot(E_\idot) = \tau_{\geq 2n-1}cp_\idot(E_\idot)^{[p]},
\qquad n \in \Z.
\end{equation}
\end{defn}

Note that since the complex $cp_\idot(E_\idot)$ is by definition
$2$-periodic, the conjugate filtration is periodic: we have
\begin{equation}\label{V.peri}
V^ncp_\idot(E_\idot) \cong V^0cp_\idot(E_\idot)[2n]
\end{equation}
for any integer $n$. By definition, the conjugate filtration is a
shift of a rescaling of the filtration $\tau^\hdot$, so that it is
commensurable to the $p$-th rescaling of the standard filtration on
$cp_\idot(E_\idot)$. This is in turn commensurable to the standard
filtration itself. Therefore by Lemma~\ref{cp.bcp.le}, the
co-periodic cyclic homology complex $\bCP_\idot(E_\idot)$ is
isomorphic to the completion of the complex $cp_\idot(E_\idot)$ with
respect to the conjugate filtration \eqref{conj.glo}. We then have a
convergent spectral sequence
\begin{equation}\label{conj.sp.glo}
H_\idot(\gr^0_Vcp_\idot(E_\idot))((u^{-1})) \Rightarrow
\bHP_\idot(E_\idot),
\end{equation}
where as before, $H_\idot(-)$ stands for homology objects, $u$ is a
formal generator of cohomological degree $2$, and we have used the
identifications \eqref{V.peri}.

\subsection{Edgewise subdivision.}

In general, the spectral sequence \eqref{conj.sp.glo} does not seem
to be particularly useful, since its initial term is rather
obscure. However, under some assumptions, it can be computed
explicitly. The first step in this computation is the so-called {\em
  edgewise subdivision}.

Recall that for any integer $l$, the category $\Lambda$ has a cousin
$\Lambda_l$ corresponding to the $l$-fold cover $S^1 \to
S^1$. Objects in $\Lambda_l$ are objects $[nl] \in \Lambda$, $n \geq
1$, equipped with the order-$l$ automorphism $\tau=\sigma^n:[nl] \to
[nl]$, and morphisms are morphisms in $\Lambda$ that commute with
the automorphism $\tau$. It is convenient to number objects in
$\Lambda_l$ by positive integers, so that $[n] \in \Lambda_l$
corresponds to $[nl] \in \Lambda$. We have the forgetful functor
$i_l:\Lambda_l \to \Lambda$, $[n] \mapsto [nl]$, and we also have a
natural projection $\pi_l:\Lambda_l \to \Lambda$ that sends $[nl]$
to the object $[n]$ considered as an induced cellular decomposition
of $S^1 \cong S^1/\tau$. The functor $\pi_l:\Lambda_l \to \Lambda$
is a Grothendieck bifibration with fiber $\ppt_l$, the groupoid with
one object with automorphism group $\Z/l\Z$. The functor $i_l$
induces the pullback functor $i_l^*:\Fun(\Lambda,R) \to
\Fun(\Lambda_l,R)$, classically known as ``edgewise subdivision
functor'', and the functor $\pi_l$ induces functors
$\pi_{l!},\pi_{l*}:\Fun(\Lambda_l,R) \to \Fun(\Lambda,R)$. We will
need a slightly more complicated functor $\pi_{l\flat}$ that sends
complexes in $\Fun(\Lambda_l,R)$ to complexes $\Fun(\Lambda,R)$.  It
is given by
\begin{equation}\label{pi.flat}
\pi_{l\flat}E_\idot = \per(\pi_{l!}(i_l^*\K_\idot \otimes E_\idot)).
\end{equation}
Equivalently, one can use $\pi_{l*}$ --- we have a natural trace map
$\pi_{l!} \to \pi_{l*}$, and when evaluated on the complex
$i_l^*\K_\idot \otimes E_\idot$, this map is an isomorphism. The
functors $i_l^*$ and $\pi_{l\flat}$ extend to filtered complexes in
the obvious way.

\medskip

Now assume that our base ring $R$ is annihilated by a prime $p$, and
restrict our attention to the $p$-fold cover $\Lambda_p$. Then we
have the following result.

\begin{prop}[{{\cite[Proposition 4.4]{ka4}}}]\label{edge.prop}
Assume given a complex $E_\idot$ in $\Fun(\Lambda,R)$ equipped with
some filtration $F^\hdot$, and denote by $E_\idot^{[p]}$ the same
complex equipped with the $p$-the rescaling $F^\hdot_{[p]}$ of the
filtration. Then there exists a functorial map
$$
\nu_p:cc_\idot(\pi_{p\flat}i_p^*E_\idot^{[p]}) \to
cp_\idot(E_\idot)^{[p]},
$$
and this map is a filtered quasiisomorphism with respect to the
standard filtrations.\endproof
\end{prop}

For the next step, we need to impose a condition on the complex
$E_\idot$.  Recall that we have assumed $pR=0$. Therefore the
cohomology algebra $H^\hdot(\Z/p\Z,R)$ of the cyclic group $\Z/p\Z$
is given by
\begin{equation}\label{c.coho}
H^\hdot(\Z/p\Z,R) \cong R[u]\langle \eps \rangle,
\end{equation}
where $u$ is a generator of degree $2$, and $\eps$ is a generator of
degree $1$. The relations are: $u$ commutes with $\eps$, and
\begin{equation}\label{p.odd}
\eps^2 = \begin{cases}
u, &\quad p=2,\\
0, &\quad p \text{ is odd.}
\end{cases}
\end{equation}
In particular, for any $R[\Z/p\Z]$-module $E$ and any integer $i$,
we have a natural map
\begin{equation}\label{eps.i}
\eps_i:\vH_i(\Z/p\Z,E) \to \vH_{i-1}(\Z/p\Z,E),
\end{equation}
where $\vH_\idot(\Z/p\Z,-)$ is the Tate homology of $\Z/p\Z$.  Since
the Tate homology of the cyclic group is $2$-periodic, we have
$\eps_i = \eps_{i+2n}$ for any $n$, so that $\eps_i$ only depends on
the parity of the integer $i$.

\begin{defn}\label{tight.def}
An $R[\Z/p\Z]$-module $E$ is {\em tight} if $\eps_1$ is an
isomorphism. A complex $E_\idot$ of $R[\Z/p\Z]$-modules is {\em
  tight} if $E_i$ is a tight $R[\Z/p\Z]$-module when $i$ is
divisible by $p$ and a projective $R[\Z/p\Z]$-module otherwise. A
complex $E_\idot$ in the category $\Fun(\Lambda_p,R)$ is {\em tight}
if for any $[n] \in \Lambda_p$, $E_\idot([n])$ is a tight complex
with respect to the action of the group $\Z/p\Z$ generated by $\tau
\in \Aut([n])$.
\end{defn}

Note that if $p$ is odd, then \eqref{p.odd} shows that for a tight
$R[\Z/p\Z]$-module $E$, we have $\eps_0=0$. Conversely, if $p=2$,
$\eps_0$ is also an isomorphism, and in fact the tightness condition
is always satisfied (both for an object and for a complex).

\begin{lemma}\label{tight.le}
\begin{enumerate}
\item Assume given a tight complex $E_\idot$ of $R[\Z/p\Z]$-modules, and
denote by $E^{[p]}_\idot$ the complex $E_\idot$ equipped with the
$p$-th rescaling of the stupid filtration. Denote by
$\vC_\idot(\Z/p\Z,E^{[p]}_\idot)$ the sum-total complex of the Tate
homology bicomplex of the group $\Z/p\Z$ with coefficients in
$E_\idot$ equipped with the filtration $F^\hdot$ induced by the filtration on
$E^{[p]}_\idot$. Let
$$
\I(E_\idot) = \tau^F_{[0,0]}\vC_\idot(\Z/p\Z,E^{[p]}_\idot).
$$
Then the induced filtration $F^\hdot$ on $\I(E_\idot)$ is the stupid
filtration, and for any integer $i$, we have a natural filtered
isomorphism
\begin{equation}\label{i.n}
\tau^F_{[i,i]}\vC_\idot(\Z/p\Z,E^{[p]}_\idot) \cong \I(E_\idot)[i].
\end{equation}
\item Assume given a tight complex $E_\idot$ in $\Fun(\Lambda_p,R)$,
  and denote by $E^{[p]}_\idot$ the complex $E_\idot$ equipped with
  the $p$-th rescaling of the stupid filtration. Consider the
  induced filtration on the complex $\pi_{p\flat}E^{[p]}_\idot$, and
  let $\I(E_\idot) = \tau^F_{[0,0]}\pi_{p\flat}E^{[p]}_\idot$. Then
  for any $[n] \in \Lambda$, we have $\I(E_\idot)([n]) \cong
  \I(E_\idot([n]))$, and for any integer $[i]$, the isomorphisms
  \eqref{i.n} induce an isomorphism
\begin{equation}\label{tau.i}
\tau^F_{[i,i]}\pi_{p\flat}E_\idot^{[p]} \cong \I(E_\idot)[i].
\end{equation}
\end{enumerate}
\end{lemma}

\proof{} Almost all of the statements are obvious; the non-obvious
ones are \cite[Lemma 5.3]{ka4}.
\endproof

Explicitly, the isomorphisms \eqref{i.n} can be described as
follows. By periodicity, we have an isomorphism
$$
u:\vC_\idot(\Z/p\Z,E^{[p]}_\idot) \cong
\vC_\idot(\Z/p\Z,E^{[p]}_\idot)[2]
$$
corresponding to the action of the generator $u$ of the cohomology
algebra \eqref{c.coho}. Twists by powers of $u$ provide isomorphisms
\eqref{i.n} for even $i$. To obtain the isomorphisms for odd $i$,
one considers the action of the generator $\eps$. This gives natural
maps 
\begin{equation}\label{eps.i.F}
\eps_i:\tau^F_{[i,i]}\vC_\idot(\Z/p\Z,E^{[p]}_\idot) \to
\tau^F_{[i-1,i-1]}\vC_\idot(\Z/p\Z,E^{[p]}_\idot),
\end{equation}
a filtered refinement of \eqref{eps.i}. Since $E_\idot$ is tight,
$\eps_i$ is a filtered isomorphism for any odd $i$.

All of this works relatively over the category $\Lambda$; in
particular, we have natural maps
\begin{equation}\label{eps.i.la}
\eps_i:\tau^F_{[i,i]}\pi_{p\flat}E^{[p]}_\idot \to
\tau^F_{[i-1,i-1]}\pi_{p\flat}E^{[p]}_\idot
\end{equation}
for any tight complex $E_\idot$ in $\Fun(\Lambda_p,R)$, and these
maps are invertible for odd $i$.

\subsection{Localized conjugate filtration.}

Now, as it turns out, for a complex $E_\idot$ in $\Fun(\Lambda,R)$
with tight edgewise subdivision $i_p^*E_\idot$, we can localize the
conjugate filtration \eqref{conj.glo} with respect to the category
$\Lambda$ and express it in terms of the complex
$\pi_{p\flat}i_p^*E_\idot$. Namely, introduce the following.

\begin{defn}\label{conj.loc.def}
For any complex $E_\idot$ in $\Fun(\Lambda_p,R)$, the {\em conjugate
  filtration} $V^\hdot$ on the complex $\pi_{p\flat}E_\idot$ is
given by
\begin{equation}\label{conj.loc}
V^n\pi_{p\flat}E_\idot = \tau^F_{\geq 2n}\pi_{p\flat}E^{[p]}_\idot,
\end{equation}
where $E^{[p]}_\idot$ stands for $E_\idot$ equipped with the $p$-the
rescaling of the stupid filtrarion.
\end{defn}

We can then take a complex $E_\idot$ in $\Fun(\Lambda,R)$, consider
the corresponding complex $\pi_{p\flat}i_p^*E_\idot$ equipped with
the conjugate filtration \eqref{conj.loc}, and apply to it the
cyclic homology complex functor $CC_\idot(-)$. Since the functor
$CC_\idot(-)$ is exact, the conjugate filtration induces a
filtration $V^\hdot$ on $CC_\idot(\pi_{p\flat}i_p^*E_\idot)$. We
denote by $\wh{CC}_\idot(\pi_{p\flat}i_p^*E_\idot)$ the completion
of the complex $CC_\idot(\pi_{p\flat}i_p^*E_\idot)$ with respect to
the filtration $V^\hdot$.

Consider now the composition
\begin{equation}\label{glo.loc}
\begin{CD}
CC_\idot(\pi_{p\flat}i_p^*E_\idot) @>{\alpha}>>
cc_\idot(\pi_{p\flat}i_p^*E_\idot) @>{\nu_p}>> cp_\idot(E_\idot)
\end{CD}
\end{equation}
of the map $\alpha$ of \eqref{CC.cc} and the natural map $\nu_p$ of
Proposition~\ref{edge.prop}.

\begin{lemma}\label{cc.v.le}
\begin{enumerate}
\item Assume that the complex $i_p^*E_\idot$ in $\Fun(\Lambda_p,R)$
  is tight in the sense of Definition~\ref{tight.def}. Then the map
  \eqref{glo.loc} extends to a quasiisomorphism
$$
\wh{CC}_\idot(\pi_{p\flat}i_p^*E_\idot) \cong \bCP_\idot(E_\idot).
$$
\item Moreover, assume that the prime $p$ is odd. Then
  \eqref{glo.loc} itself is a filtered quasiisomorphism with respect
  to the filtrations $V^\hdot$ of \eqref{conj.loc}
  resp.\ \eqref{conj.glo}.
\item In addition, still assuming that $p$ is odd, let
  $\I(i_p^*E_\idot)$ be canonical complex of
  Lemma~\ref{tight.le}. Then we have a natural isomorphism
\begin{equation}\label{gr.v}
\gr^0_V\pi_{p\flat}i_p^*E_\idot \cong \K(\I(i_p^*E_\idot))
\end{equation}
in the derived category $\D(\Lambda,R)$.
\end{enumerate}
\end{lemma}

\proof{} \thetag{i} is \cite[Lemma 5.19]{ka4}, \thetag{ii} is
\cite[Lemma 5.9]{ka4}, and \thetag{iii} is \cite[Lemma
  5.7]{ka4}.
\endproof

By virtue of this result, for an odd prime $p$, we can rewrite the
spectral sequence \eqref{conj.sp.glo} as
\begin{equation}\label{conj.sp.loc}
CC_\idot(\gr^0_V\pi_{p\flat}i_p^*E_\idot)((u^{-1})) \Rightarrow
\bHP_\idot(E_\idot),
\end{equation}
where we have used the obvious counterpart of the periodicity
isomorphism \eqref{V.peri} for the filtration
\eqref{conj.loc}. Moreover, by \eqref{hc.hla}, \eqref{gr.v}, and
\eqref{hc.hh}, we can further rewrite \eqref{conj.sp.i} as
\begin{equation}\label{conj.sp.i}
HH_\idot(\I(i_p^*E_\idot))((u^{-1})) \Rightarrow
\bHP_\idot(E_\idot),
\end{equation}
where as in \eqref{conj.sp.glo}, $u$ is a formal generator of
cohomological degree $2$.

\subsection{DG algebras.}

To make the spectral sequence \eqref{conj.sp.i} useful, it remains
to compute the complex $\I(i_p^*E_\idot)$. In order to do this, we
need to assume further that $E_\idot$ comes from a DG algebra
$A_\idot$. We thus assume given a commutative ring $k$ annihilated
by an odd prime $p$. We denote by $k^{(1)}$ the Frobenius twist of
$k$ --- that is, $k$ considered as a module over itself via the
absolute Frobenius map $k \to k$. For any flat $k$-module $V$, we
denote $V^{(1)} = V \otimes_k k^{(1)}$.

\begin{prop}[{{\cite[Proposition 6.10]{ka4}}}]\label{ten.prop}
Assume that the ring $k$ is Noetherian. Then for any complex
$V_\idot$ of flat $k$-modules, the complex $V_\idot^{\otimes_k p}$
of $k[\Z/p\Z]$-modules is tight in the sense of
Definition~\ref{tight.def}, and we have a natural identification
$$
\I(V_\idot^{\otimes_k p}) \cong V_\idot^{(1)},
$$
where $\I(-)$ is the canonical complex provided by
Lemma~\ref{tight.le}.\endproof
\end{prop}

\begin{corr}\label{ten.corr}
Assume given a DG algebra $A_\idot$ termwise-flat over a commutative
ring $k$ annihilated by an odd prime $p$. Then the complex
$i_p^*A_\idot^\hush$ in the category $\Fun(\Lambda_p,k)$ is tight in
the sense of Definition~\ref{tight.def}, and we have a natural
identification
$$
\I(i_p^*A_\idot^\hush) \cong \left(A_\idot^\hush\right)^{(1)}.
$$
\end{corr}

\proof{} This is \cite[Lemma 6.19]{ka4}.
\endproof

By virtue of this corollary, we have a natural identification
$$
HH_\idot(\I(i_p^*A^\hush_\idot)) \cong HH^{(1)}_\idot(A_\idot),
$$
where for any complex $E_\idot$ in $\Fun(\Lambda,k)$ termwise-flat
over $k$, $HH^{(1)}_\idot(E_\idot)$ denotes the homology of the
Frobenius twist $CH_\idot(E_\idot)^{(1)}$ of the Hochschild homology
complex $CH_\idot(E_\idot)$. If the commutative ring $k$ is finitely
generated and regular, so that $k^{(1)}$ is flat over $k$, then we
have
$$
HH^{(1)}_\idot(A_\idot) \cong HH_\idot(A^{(1)}_\idot),
$$
where $A^{(1)}_\idot$ is the Frobenius twist of the DG algebra
$A_\idot$. Then \eqref{conj.sp.i} takes its final form --- what we
have is a spectral sequence
\begin{equation}\label{conj.sp}
HH_\idot(A^{(1)}_\idot)((u^{-1})) \Rightarrow \bHP_\idot(A_\idot).
\end{equation}
This is the {\em conjugate spectral sequence} for the DG algebra
$A_\idot$.

\medskip

We finish the section with an alternative description of the complex
$i_p^*A^\hush_\idot$ in terms of the $p$-tensor power algebra
$A^{\otimes_k p}_\idot$; this goes back to \cite[Subsection
  2.2]{ka0}, and we will need it later in Section~\ref{tate.sec}.

For any small category $\C$, sending an object $[n] \in \Lambda$ to
$\C^{V([n])}$ defines a contravariant functor from $\Lambda$ to the
category of small categories, and the Grothendieck construction
\cite{sga} associates a category fibered over $\Lambda$ to this
functor. We will denote the category by $\C\wr\Lambda$, and we will
denote the fibration by $\pi:\C\wr\Lambda \to \Lambda$. Then for any
DG algebra $A_\idot$ in the category $\Fun(\C,k)$ --- or
equivalently, for any functor from $\C$ to DG algebras over $k$ ---
the construction of the complex $A^\hush_\idot$ of
Subsection~\ref{cycl.subs} admits a straightforward refinement that
produces a complex $A^\hush_\idot$ in the category
$\Fun(\C\wr\Lambda,k)$.

Consider the case $\C=\ppt_p$, the groupoid with one object with
automorphism group $G=\Z/p\Z$. Then a DG algebra in $\Fun(\ppt_p,k)$
is simply a DG algebra $A_\idot$ over $k$ equipped with an action of
the group $G$, and what we obtain is a complex $A^\hush_\idot$ in
the category $\Fun(\ppt_p\wr\Lambda,k)$. The categories $\Lambda$
and $\ppt_p\wr\Lambda$ have the same objects, and for any such
object $[n]$, we have
$$
A^\hush_\idot([n])=A^{\otimes_k n}_\idot,
$$
with the natural action of $G^n$. By base change \cite[Lemma
  1.7]{ka0}, we have
\begin{equation}\label{hush.n}
\pi_*A^\hush_\idot([n]) \cong \left(A^\hush_\idot([n])\right)^{G^n} =
\overline{A}^\hush_\idot([n]),
\end{equation}
where we denote by $\overline{A}_\idot = A^G_\idot \subset A_\idot$
the subalgebra of $G$-invariants. Taken together, these
identifications produce an isomorphism
\begin{equation}\label{hush.tot}
\pi_*A^\hush_\idot \cong \overline{A}^\hush_\idot.
\end{equation}
Now note that we also have a natural embedding
\begin{equation}\label{lambda.eq}
\lambda:\Lambda_p \to \ppt_p\wr\Lambda
\end{equation}
that induces the diagonal embedding $\ppt_p \to \ppt_p^n$ on the
fiber over any object $[n] \in \Lambda$. Then for any DG algebra
$A_\idot$ over $k$, we have natural isomorphism
\begin{equation}\label{la.i}
i_p^*A_\idot^\hush \cong \lambda^*\left(A^{\otimes_k
  p}_\idot\right)^\hush,
\end{equation}
where in the right-hand side, the $p$-th power $A^{\otimes_k
  p}_\idot$ is equipped with the natural $G$-action by the longest
permutation.

\section{Tate cohomology.}\label{tate.sec}

To analyze the conjugate spectral sequence \eqref{conj.sp}, we need
a more invariant definition of the functor $\pi_{p\flat}$ of
\eqref{pi.flat}. The relevant conceptual formalism is that of Tate
cohomology.

\subsection{Relative Tate cohomology.}

Assume given a ring $k$, a finite group $G$, and a bounded complex
$E_\idot$ of $k[G]$-modules. Recall that the {\em Tate cohomology}
of the group $G$ with coefficients in $E_\idot$ is given by
\begin{equation}\label{vh.eq}
\vH^\hdot(G,E_\idot) =
\Ext^\hdot_{\D(k[G])/\D^{pf}(k[G])}(k,E_\idot),
\end{equation}
where $k$ is the trivial $k[G]$-module, and $\D(k[G])/\D^{pf}(k[G])$
is the quotient of the derived category $\D(k[G])$ of all
$k[G]$-modules by its full subcategory $\D^{pf}(k[G]) \subset
\D(k[G])$ spanned by perfect complexes of $k[G]$-modules. In order
to compute Tate cohomology, it is convenient to introduce the
following (we use the same notation and terminology as in
\cite[Subsection 6.3]{ka5}).

\begin{defn}\label{reso.def}
\begin{enumerate}
\item {\em Resolution data} for a finite group $G$ is a pair $\nu =
  \langle P_\idot,I^\hdot \rangle$ of a left free resolution
  $P_\idot$ and a right free resolution $I^\hdot$ of the trivial
  $\Z[G]$-module $\Z$.
\item For any associative unital ring $k$, any bounded complex
  $E_\idot$ of $k[G]$-modules, and any resolution data $\nu$, the
  {\em Tate cohomology complex} of $G$ with coefficients in
  $E_\idot$ is given by
\begin{equation}\label{c.phi.tate}
\vC^\hdot(G,\nu,E_\idot) = \left(E_\idot \otimes \overline{P}_\idot
\otimes I^\hdot\right)^G,
\end{equation}
where $\overline{P}_\idot$ is the cone of the augmentation map $P_\idot
\to \Z$.
\item For any associative unital ring $k$, any bounded complex
  $E_\idot$ of $k[G]$-modules, and any resolution data $\nu$, the
  {\em reduced Tate cohomology complex} of $G$ with coefficients in
  $E_\idot$ is given by
\begin{equation}\label{c.phi.tate.bis}
\vC^\hdot_{red}(G,\nu,E_\idot) = \left(E_\idot \otimes
\wt{P}_\idot\right)^G,
\end{equation}
where $\wt{P}_\idot$ is the cone of the natural map $P_\idot \to \Z
\to I^\hdot$.
\end{enumerate}
\end{defn}

\begin{remark}
Here as well as elsewhere, we use homological and cohomological
indexes for complexes interchangeably, with the convention that $E_i
= E^{-i}$.
\end{remark}

Resolution data form a category in the obvious way, and this
category is connected (for example, every $\nu$ is obviously
connected by a chain of maps to a functorial resolution data set
$\nu'$ obtained by taking bar resolutions). For any resolution data
$\nu$, we have a natural map $\wt{P}_\idot \to \overline{P}_\idot
\otimes I^\hdot$, and the induced map
$$
\vC^\hdot_{red}(G,\nu,E_\idot) \to \vC^\idot(G,\nu,E_\idot)
$$
of Tate cohomology complexes is a quasiisomorphism for any bounded
complex $E_\idot$. Moreover, a map of resolution data induces a map
of Tate complexes, and these maps are also
quasiisomorphisms. Therefore we can drop $\nu$ from notation and
obtain a well-defined object $\vC^\idot(G,E_\idot) \cong
\vC^\idot_{red}(G,E_\idot)$ in the derived category of $k$-modules. Its
cohomology modules are then canonically identified with Tate
cohomology groups $\vH^\hdot(G,E_\idot)$.

\medskip

It is obvious from \eqref{vh.eq} that $\vH^\hdot(G,k)$ is an
algebra, and for any bounded complex $E_\idot$,
$\vH^\hdot(G,E_\idot)$ is a module over $\vH^\hdot(G,k)$. To lift it
to the dervied category level, one chooses resolution data
$\nu=\langle P_\idot,I^\hdot \rangle$ that are {\em multiplicative}
in the following sense: both $I^\hdot$ and $\overline{P}_\idot$ are
DG algebras, and the natural maps $\Z \to I^\hdot$, $\Z \to
\overline{P}_\idot$ are algebra maps. Multiplicative resolution data
exist. For each multiplicative $\nu$, $\vC^\hdot(G,\nu,k)$ is a DG
algebra, and $\vC^\hdot(G,\nu,E_\idot)$ is a module over this DG
algebra. Passing to the derived category, we obtain well-defined
multiplication maps
$$
\vC^\hdot(G,k) \lotimes \vC^\hdot(G,E_\idot) \to
\vC^\hdot(G,E_\idot)
$$
that do not depend on the choice of the resolution data $\nu$.

\medskip

Assume now given small categories $\C$, $\C'$ and a functor $\pi:\C'
\to \C$. Assume further that $\pi$ is a Grothendieck bifibration in
the sense of \cite{sga}, and that the fibers of this bifibrations
are equivalent to $\ppt_G$, the groupoid with one object with
automorphism group $G$. By base change \cite[Lemma 1.7]{ka0}, for
any object $c' \in \C'$ with $c=\pi(c') \in \C$, and any functor $E
\in \Fun(\C',k)$, we have natural isomorphisms
\begin{equation}\label{bc}
\pi_!E(c') \cong E(c')_G, \qquad \pi_*E(c) \cong E(c')^G,
\end{equation}
where $G$ acts on $E(c')$ via the embedding $\ppt_G \to \C'$ of the
fiber over $c \in \C$. Then one can obviously make the constructions
above work ``relatively over $\C$''. Namely, one defines {\em resolution
  data} for $\pi$ as pairs $\nu = \langle P_\idot,I^\hdot \rangle$
of a left and a right resolution of the constant functor $\Z \in
\Fun(\C',\Z)$ such that for any $c \in \C$, $\nu$ restricted to the
fiber $\ppt_G \subset \C'$ over $c$ gives resolution data in the
sense of Definition~\ref{reso.def}. One shows easily that resolution
data exist (e.g. take the bar resolutions) and that the category of
resolution data is connected. Then for any resolution data $\nu$,
one defines
\begin{equation}\label{pi.nu.fl}
\pi_\nu^\flat(E) = \pi_*(E \otimes \overline{P}_\idot \otimes
I^\hdot), \qquad \pi_\nu^\flfl(E) = \pi_*(E \otimes
\wt{P}_\idot),
\end{equation}
where $\overline{P}_\idot$ is as in \eqref{c.phi.tate}, and
$\wt{P}_\idot$ is as in \eqref{c.phi.tate.bis}. By \eqref{bc}, we
have a natural quasiisomorphism $\pi_\nu^\flat(E) \cong
\pi_\nu^\flfl(E)$, and both complexes do not depend on the choice
of $\nu$. All in all, we obtain a well-defined object
$$
\pi^\flat(E) \in \D(\C,k)
$$
in the derived category $\D(\C,k)$. For any object $c' \in \C'$ with
$c = \pi(c') \in \C$, \eqref{bc} gives a natural identification
$$
\pi^\flat(E)(c) \cong \vC^\hdot(G,E(c')).
$$
Moreover, as in \cite[Subsection 6.3]{ka5}, one can choose
multiplicative resolution data for $\pi$. This shows that
$\pi^\flat(k)$ is an algebra object in $\D(\C,k)$, and for any $E$,
we have natural action maps
\begin{equation}\label{mult.gen}
\pi^\flat(k) \lotimes \pi^\flat(E) \to \pi^\flat(E)
\end{equation}
that turn $\pi^\flat(E)$ into a module object over the algebra
$\pi^\flat(k)$.

Finally, we observe that all of the above can be repeated {\em
  verbatim} for a bounded complex $E_\idot$ instead of a single
object $E$ in $\Fun(\C',k)$. Moreover, assume given a filtered
complex $\langle E_\idot,F^\hdot \rangle$ in $\Fun(\C',k)$, and
assume that the filtration $F^\hdot$ is termwise-split, and
$\gr^i_FE_\idot$ is bounded for any integer $i$. Then again,
\eqref{pi.nu.fl} gives an object
\begin{equation}\label{pi.fl.gen}
\pi^\flat E_\idot \in \DF(\C,k)
\end{equation}
in the filtered derived category $\DF(\C,k)$, and as such, it does
not depend on the choice of the resolution data $\nu$.

\subsection{Conjugate filtration.}

We now want to apply relative Tate cohomology to the study of the
conjugate filtration of Definition~\ref{c.phi.tate}. We assume that
the base ring $k$ is commutative and annihilated by a prime $p$. For
simplicity, we also assume right away that $k$ is finitely generated
and regular, so that the absolute Frobenius map $k \to k$ is finite
and flat. Now, consider the natural projection
$$
\pi_p:\Lambda_p \to \Lambda.
$$
This is a Grothendieck bifibration, and its fiber is $\ppt_p$, the
groupoid with one object with automorphism group $\Z/p\Z$. Therefore
for any filtered complex $\langle E_\idot,F^\hdot \rangle$ in
$\Fun(\Lambda_p,k)$ such that the filtration $F^\hdot$ is
termwise-split and $\gr^i_FE_\idot$ is bounded for any $i$,
\eqref{pi.fl.gen} provides a natural object
\begin{equation}\label{pi.p.f}
\pi_p^\flat E_\idot \in \DF(\Lambda,k).
\end{equation}
We note that we have a natural isomorphism
\begin{equation}\label{pi.pi}
\pi_\nu^\flfl E_\idot \cong \pi_{p\flat}E_\idot[1],
\end{equation}
where $\pi_{p\flat}E_\idot$ is the complex of \eqref{pi.flat}, and
$\nu$ is the pair of resolutions of $\Z$ obtained by periodization
of the complex $\K_\idot$ of \eqref{4term}. Therefore the object
\eqref{pi.p.f} coincides with $\pi_{p\flat}E_\idot$ up to a
homological shift.

In particular, let $E_\idot$ be an arbitrary complex in
$\Fun(\Lambda_p,k)$, and let $E^{[p]}_\idot$ be $E_\idot$ equipped
with the $p$-th rescaling of the stupid filtration. Then all the
assumptions on the filtration are satisfied, so that we have a
well-defined object $\pi_p^\flat E^{[p]}_\idot$ in
$\DF(\Lambda,k)$. We also have the action map
\begin{equation}\label{mult.pi}
\pi_p^\flat k \lotimes \pi_p^\flat E^{[p]}_\idot \to
\pi_p^\flat E^{[p]}_\idot
\end{equation}
induced by the map \eqref{mult.gen}.

\begin{lemma}[{{\cite[Lemma 3.2]{ka-trudy}}}]
We have a canonical isomorphism
\begin{equation}\label{sum.pi}
\pi_{p\flat} k \cong \K_\idot(k)((u)) = \bigoplus_{i \in
  \Z}\K_\idot(k)[2i]
\end{equation}
of complexes in $\Fun(\Lambda,k)$.\endproof
\end{lemma}

By \eqref{pi.pi}, this gives an isomorphism $\tau_{[i,i]}\pi_p^\flat
k \cong k[i]$ for any integer $i$. Since the filtered truncation
functors \eqref{tau.F} are obviously compatible with the tensor
products, the map \eqref{mult.pi} then induces a map
$$
\eps_i:\tau_{[1,1]}\pi_p^\flat k \otimes
\tau^F_{[i,i]}\pi_p^\flat E^{[p]}_\idot \cong 
\tau^F_{[i,i]}\pi_p^\flat E^{[p]}_\idot[1] \to
\tau^F_{[i+1,i+1]}\pi_p^\flat E^{[p]}_\idot
$$
for any integer $i$. These are exactly the maps \eqref{eps.i.la}. In
particular, $\eps_i$ only depends on the parity of $i$, and if
$E_\idot$ is tight in the sense of Definition~\ref{tight.def}, then
$\eps_1$ is an isomorphism. In this case, for any integer $j$,
\eqref{mult.pi} also induces an isomorphism
\begin{equation}\label{i.n.mult}
\I(E_\idot)[j] = \tau_{[j,j]}\pi^p_\flat k \otimes
\tau^F_{[0,0]}\pi_p^\flat E^{[p]}_\idot \cong
\tau^F_{[j,j]}\pi_p^\flat E_\idot^{[p]},
\end{equation}
a version of the isomorphism \eqref{i.n}.

If $p$ is odd, tightness of $E_\idot$ further implies that
$\eps_0=0$. One immediate corollary of this is a short construction
of the quasiisomorphism \eqref{gr.v} of Lemma~\ref{cc.v.le}, for odd
$p$. Indeed, since $\eps_2=\eps_0=0$ and the filtered truncations
are multiplicative, the map \eqref{mult.pi} induces a natural map
$$
\begin{aligned}
\K(\I(E_\idot))[1] &= \K_\idot[1] \otimes \I(E_\idot) \cong
\tau_{[1,2]}\pi_p^\flat k \otimes
\tau^F_{[0,0]}\pi_p^\flat E^{[p]}_\idot \to\\
&\to \tau^F_{[1,2]}\pi_p^\flat E^{[p]}_\flat \cong
\gr^0_V\pi_{p\flat}E_\idot[1],
\end{aligned}
$$
and since the maps \eqref{i.n.mult} are isomorphisms, this map is a
quasiisomorphism.

However, we will need another corollary. Namely, keep the assumption
that $E_\idot$ is tight, let
\begin{equation}\label{pi.bar}
\overline{\pi}^\flat_p E_\idot = \tau^F_{\geq 0}\pi_p^\flat
E^{[p]}_\idot,
\end{equation}
and consider this complex as an object in the filtered derived
category $\DF(\Lambda,k)$ by equipping it with the filtration
$\tau^\hdot$, via the functor \eqref{F.tau}. We then have the
augmentation map
\begin{equation}\label{I.aug}
a:\overline{\pi}_p^\flat E_\idot \to
\tau_{[0,0]}\pi_p^\flat E^{[p]}_\idot = \I(E_\idot)
\end{equation}
in $\DF(\Lambda,k)$, where $\I(E_\idot)$ is placed in filtered
degree $0$.

\begin{lemma}\label{spl.dege.le}
Assume that the map $a$ of \eqref{I.aug} admits a one-sided inverse
$s:\I(E_\idot) \to \overline{\pi}_p^\flat E_\idot$, $s
\circ a = \id$ in the filtered derived category
$\DF(\Lambda,k)$. Then the spectral sequence \eqref{conj.sp.i}
degenerates.
\end{lemma}

\proof{} By definition, we have a natural embedding
$\overline{\pi}_p^\flat E_\idot \to \pi_p^\flat E_\idot$. Composing
it with $s$, we obtain a map
$$
\I(E_\idot) \to \pi_p^\flat E_\idot,
$$
and by \eqref{mult.pi}, this map induces a map
$$
b:\pi_p^\flat k \otimes \I(E_\idot) \cong \pi_p^\flat E_\idot.
$$
For any integer $i$, the associated graded quotient $\gr^i_\tau(b)$
is the isomorphism \eqref{tau.i} of Lemma~\ref{tight.le}, so that
$b$ is a filtered quasiisomorphism. By the direct sum decomposition
\eqref{sum.pi}, this implies that the conjugate filtration on
$\pi^p_\flat E_\idot$ splits.
\endproof

\subsection{Splitting for DG algebras.}

Now keep the assumptions of the previous subsection, and assume
given a DG algebra $A_\idot$ termwise-flat over the commutative ring
$k$. Recall that by Corollary~\ref{ten.corr}, the corresponding
complex $i_p^*A^\hush_\idot$ is tight in the sense of
Definition~\ref{tight.def}, so that Lemma~\ref{spl.dege.le}
applies. To finish the section, we prove one corollary of this fact.

Choose resolution data $\nu=\langle P_\idot,I^\hdot \rangle$ for the
group $G=\Z/p\Z$ that are multiplicative, so that for any DG algebra
$B_\idot$ termwise-flat over $k$ and equipped with a $G$-action, the
Tate cohomology complex $\vC_\idot(G,B_\idot)$ is a DG algebra over
$k$. Consider the $p$-fold tensor product $A_\idot^{\otimes_k p}$,
and let $\Z/p\Z$ act on it by the longest permutation. Moreover,
equip this tensor product with the $p$-the rescaling of the stupid
filtration.

\begin{defn}
The DG algebra $P_\idot(A_\idot)$ is given by
$$
P_\idot(A_\idot) = \tau^F_{\geq 0}\vC^\hdot(\Z/p\Z,\nu,A^{\otimes_k
  p}_\idot),
$$
where $A^{\otimes_k p}_\idot$ is equipped with the $p$-th rescaling
of the stupid filtration.
\end{defn}

As for the complex \eqref{pi.bar}, we equip $P_\idot(A_\idot)$ with
the filtration $\tau^\hdot$ and treat it as a filtered DG algebra.
Since the filtered truncation functors $\tau^F_{\geq \idot}$ are
multiplicative, $P_\idot(A_\idot)$ is well-defined, and up to a
filtered quasiisomorphism, it does not depend on the choice of
resolution data $\nu$. By Proposition~\ref{ten.prop}, we have
\begin{equation}\label{gr.P}
\gr^i_\tau P_\idot(A_\idot) \cong A_\idot^{(1)}[i]
\end{equation}
for any integer $i \geq 0$. In particular, we have an augmentation
map
$$
a:P_\idot(A_\idot) \to A_\idot^{(1)},
$$
and it is a filtered DG algebra map (where $A^{(1)}_\idot$ is in
filtered degree $0$).

\begin{prop}\label{P.dege.prop}
Assume that the prime $p$ is odd, and that there exists a filtered
DG algebra $A'_\idot$ over $k$ and a filtered DG algebra map
$s:A'_\idot \to A_\idot^{(1)}$ such that the composition $s \circ
a:A'_\idot \to A^{(1)}_\idot$ is a filtered quasiisomorphism. Then
the spectral sequence \eqref{conj.sp} for the DG algebra $A_\idot$
degenerates.
\end{prop}

\proof{} Fix multiplicative resolution data $\langle P_\idot,I^\hdot
\rangle$ for the group $G=\Z/p\Z$, and consider the cone
$\overline{P}_\idot$ of the augmentation map $P_\idot \to \Z$. By
definition, both $I^\hdot$ and $\overline{P}_\idot$ are DG algebras
over $k$ equipped with a $G$-action.

For every integer $n \geq 1$, the complex
$\overline{P}_\idot^{\otimes n}$ is concentrated in non-negative
homological degrees, its degree-$0$ term is $\Z$, while all the
other terms are free $\Z[G]$-modules. Therefore
$\overline{P}_\idot^{\otimes n} = \overline{P}^n_\idot$ for some
free left resolution $P^n_\idot$ of the trivial module $\Z$. The
complex $I^{\hdot\otimes n}$ is a right free resolution of $\Z$, so
that $\langle P^n_\idot,I^{\hdot\otimes n}\rangle$ also gives
resolution data for $G$ in the sense of
Definition~\ref{reso.def}.

Now, since both $\overline{P}_\idot$ and $I^\hdot$ are DG algebras
equipped with a $G$-action, we have natural complexes
$\overline{P}^\hush_\idot$, $(I^\hdot)^\hush$ in the category
$\Fun(\ppt_p\wr\Lambda,k)$. Restricting them to $\Lambda_p \subset
\ppt_p\wr\Lambda$ with respect to the embedding \eqref{lambda.eq},
we obtain complexes
$$
\overline{P}^\lambda_\idot = \lambda^*\overline{P}_\idot^\hush,
\qquad I^\hdot_\lambda = \lambda^*(I^\hdot)^\hush
$$
in the category $\Fun(\Lambda_p,k)$, and for any object $[n] \in
\Lambda_p$, we have natural identifications
$$
\overline{P}^\lambda_\idot([n]) \cong \overline{P}_\idot^{\otimes
  n}, \qquad I^\hdot_\lambda([n]) \cong I^{\hdot\otimes n}.
$$
Thus we can put together resolution data $\langle
P^n_\idot,I^{\hdot\otimes n} \rangle$, $n \geq 1$ for the group $G$
into resolution data $\langle P^\lambda_\idot,I^\hdot_\lambda
\rangle$ for the bifibration $\pi_p:\Lambda_p \to \Lambda$ such that
$\overline{P^\lambda}_\idot \cong \overline{P}^\lambda_\idot$, and
these resolution data can be then used for computing the relative
Tate cohomology functor $\pi_p^\flat$. This gives a natural
identification
\begin{equation}\label{pi.p.1}
\pi_p^\flat i_p^*A^\hush \cong \pi_{p*}\left(i_p^*A^\hush_\idot
\otimes \overline{P}^\lambda_\idot \otimes I^\hdot_\lambda\right).
\end{equation}
Now denote $B_\idot = A^{\otimes_k p}_\idot \otimes
\overline{P}_\idot \otimes I^\hdot$, and consider it as a
$G$-equivariant DG algebra over $k$. Then by virtue of \eqref{la.i},
we can rewrite \eqref{pi.p.1} as
\begin{equation}\label{pi.p.2}
\pi_p^\flat i_p^*A^\hush_\idot \cong \pi_{p*}\lambda^*B^\hush_\idot.
\end{equation}
Note that we have a natural map
\begin{equation}\label{la.pi}
\pi_*B^\hush_\idot \to \pi_{p*}\lambda^*B^\hush_\idot,
\end{equation}
where $\pi:\ppt_p\wr\Lambda \to \Lambda$ is the natural fibration.
At each object $[n] \in \Lambda$, we can evaluate
$\pi_*B^\hush_\idot$ by \eqref{hush.n}, and then this map is just
the natural embedding
$$
\left(B^G_\idot\right)^{\otimes_k n} = \left(B^{\otimes_k
  n}_\idot\right)^{G^n} \to \left(B^{\otimes_k n}_\idot\right)^G.
$$
Moreover, by \eqref{hush.tot} and \eqref{pi.p.2}, the map
\eqref{la.pi} actually gives a natural map
\begin{equation}\label{la.pi.bis}
\overline{B}_\idot^\hush \cong \pi_*B^\hush_\idot \to
\pi_p^\flat i_p^*A^\hush_\idot.
\end{equation}
where $\overline{B}_\idot \subset B_\idot$ is the subalgebra of
$G$-invariants.

Now equip $A^{\otimes_k p}_\idot$ with the $p$-the rescaling of the
stupid filtration, and consider the corresponding filtrations on the
algebras $B_\idot$, $\overline{B}_\idot$. Then since $\tau^F_{\geq
  0}$ is a multiplicative functor, the natural map \eqref{la.pi.bis}
induces a map
\begin{equation}\label{la.pi.3}
(\tau^F_{\geq 0}\overline{B}_\idot)^\hush \to \tau^F_{\geq
    0}\pi_p^\flat i_p^*A^\hush_\idot = \overline{\pi}^p_\flat
  i_p^*A^\hush_\idot.
\end{equation}
By construction, if we compose this map with the projection
$$
\overline{\pi}_p^\flat i_p^*A^\hush_\idot \to
\gr^0_\tau\overline{\pi}_p^\flat i_p^*A^\hush_\idot \cong
\I(i_p^*A^\hush_\idot) \cong A^{(1)\hush}_\idot,
$$
then the resulting map is induced by the augmentation map
$$
\tau^F_{\geq 0}\overline{B}_\idot \to
\tau^F_{[0,0]}\overline{B}_\idot \cong A^{(1)}_\idot.
$$
It remains to spell out the notation. By definition, we actually
have
$$
\overline{B}_\idot = \vC_\idot(G,\nu,A^{\otimes_k n}),
$$
where $\nu$ stands for our original resolution data $\langle
P_\idot,I^\hdot \rangle$. Therefore $\tau^F_{\geq
  0}(\overline{B}_\idot) \cong P_\idot(A_\idot)$, and the map
\eqref{la.pi.3} is a map
$$
P_\idot(A_\idot)^\hush \to A^{(1)}_\hush.
$$
If there exists a filtered DG algebra $A'_\idot$ and a filtered DG
algebra map $s:A'_\idot \to P_\idot(A_\idot)$ satisfying the
assumptions of the Proposition, then we also have a filtered map
$A^{'\hush}_\idot \to P_\idot(A_\idot)^\hush$, and the composition
map
$$
\begin{CD}
A^{'\hush}_\idot @>>> P_\idot(A_\idot) @>>> \overline{\pi}_p^\flat
i_p^*A^\hush_\idot @>>> A^{(1)\hush}_\idot
\end{CD}
$$
is a filtered quasiisomorphism. Then we are done by
Lemma~\ref{spl.dege.le}.
\endproof

\section{Characteristic $2$.}\label{2.sec}

We now make a digression and explain how to modify the arguments of
\cite{ka4} to obtain the conjugate spectral sequence \eqref{conj.sp}
in characteristic $2$. The problem here is Lemma~\ref{cc.v.le}:
while \thetag{i} is true in any characteristic, \thetag{ii} is
definitely wrong in characteristic $2$, and it is currently unknown
whether \thetag{iii} is true or not. However, there is the following
weaker result.

\begin{prop}\label{char.2.prop}
For any DG algebra $A_\idot$ over a perfect field $k$ of positive
characteristic $p$, there exists a natural isomorphism
\begin{equation}\label{char.2}
HC_\idot(\gr^0_V\pi_{p\flat}i_p^*A^\hush_\idot) \cong
HH_\idot(A^{(1)}_\idot).
\end{equation}
\end{prop}

While weaker than \eqref{gr.v}, this identification still allows one
to rewrite \eqref{conj.sp.loc} in the form \eqref{conj.sp}, at least
for DG algebras over a perfect field. For degeneration questions,
this is irrelevant; the reader who is only integerested in
degeneration of the spectral sequences can safely skip this section.

\subsection{Trace functors.}

To get a better handle on the complex
$\pi_{p\flat}i_p^*A^\hush_\idot$, we use the formalism of trace
functors of \cite[Section 2]{ka6}. Here are the basic ingredients.

One starts by ``categorifying'' the construction of the object
$A^\hush$ of Subsection~\ref{cycl.subs}. To every small monoidal
category $\C$, one associates a covariant functor from $\Lambda$ to
the category of small categories that sends $[n] \in \Lambda$ to
$\C^{V([n])}$, and sends a morphism $f:[n'] \to [n]$ to the product
of multiplication functors $m_{f^{-1}(v)}:\C^{f^{-1}(v)} \to \C$, $v
\in V([n])$ induced by the monoidal structure on $\C$, as in
\eqref{a.hush}. Applying the Grothendieck construction, one obtains
a category $\C^\hush$ and a cofibration $\rho:\C^\hush \to
\Lambda$. This is somewhat similar to the wreath product
construction $\C\wr\Lambda$, except that the functor is covariant,
not contravariant, and the projection $\rho:\C^\hush \to \Lambda$ is
a cofibration, not a fibration.

Explicitly, objects of $\C^\hush$ are pairs $\langle [n],\{c_v\}
\rangle$ of an object $[n] \in \Lambda$ and a collection $\{c_v\}$
of objects in $\C$ numbered by vertices $v \in V([n])$. A morphism
from $\langle [n'],\{c'_v\} \rangle$ to $\langle [n],\{c_v\}
\rangle$ is given by a morphism $f:[n'] \to [n]$ and a collection of
morphisms
$$
f_v:\bigotimes_{v' \in f^{-1}(v)}c'_{v'} \to c_v, \qquad v \in
V([n]).
$$
A morphism is {\em cartesian} if all the components $f_v$ are
invertible. Note that stated in this way, the definition makes
perfect sense even when the category $\C$ is not small.

\begin{defn}
A {\em trace functor} from a monoidal category $\C$ to some category
$\E$ is a functor $F:\C^\hush \to \E$ that sends cartesian maps in
$\C^\hush$ to invertible maps in $\E$.
\end{defn}

Explicitly, a trace functor is given by a functor $F:\C \to \E$ and
a collection of isomorphisms
\begin{equation}\label{tau.mn}
\tau_{M,N}:F(M \otimes N) \cong F(N \otimes M), \qquad M,N \in \C
\end{equation}
satisfying some compatibility constraints (see \cite[Subsection
  2.1]{ka6}). A trivial example of a trace functor is obtained by
fixing a commutative ring $k$, and taking $\C=\E=k\amod$, the
category of of $k$-modules, with $\tau_{-,-}$ being the standard
commutativity isomorphisms. There are also non-trivial examples. One
such was considered in \cite{ka6} in detail. We still take
$\C=\E=k\amod$, fix an integer $l \geq 1$, and let
\begin{equation}\label{F.cycl}
F(V) = V^{\otimes_k l}_\sigma,
\end{equation}
where $\sigma:V^{\otimes_k l} \to V^{\otimes_k l}$ is the order-$l$
permutation. The maps \eqref{tau.mn} are given by $\tau_{M,N} =
\tau'_{M,N} \circ (\sigma_M \otimes \id)$, where $\tau'_{M,N}$ are
the commutativity maps, and $\sigma_M$ is the order-$n$ permutation
acting on $M^{\otimes_k l}$.

Every algebra object $A$ in the monoidal category $\C$ defines a
section $\alpha:\Lambda \to \C^\hush$ of the cofibration
$\rho:\C^\hush \to \Lambda$, and composing this section with a trace
functor $F$ gives a natural functor
$$
FA^\hush = F \circ \alpha:\Lambda \to \E.
$$
If $\C = \E = k\amod$, what we obtain is an object $FA^\hush$ in the
category $\Fun(\Lambda,k)$ associated to any associative unital
algebra $A$ over $k$. If $F$ is the identity functor with the
tautological teace functor structure, then this just the object
$A^\hush$ of Subsection~\ref{cycl.subs}. In general, we obtain a
version of cyclic homology twisted by the trace functor $F$, the
main object of study in \cite{ka6}.

Another way to express this is to say that a trace functor $F:\C \to
k\amod$ defines an object $F^\hush \in \Fun(\C^\hush,k)$, and we
have $FA^\hush = \alpha^*F^\hush$. Analogously, a trace functor from
$\C$ to the category $C_\idot(k)$ of complexes of $k$-modules gives
a complex $F^\hush$ in $\Fun(\C^\hush,k)$, and $FA^\hush =
\alpha^*F^\hush$ is a complex in $\Fun(\Lambda,k)$.

\begin{remark}
Strictly speaking, when $\C$ is not small, $\Fun(\C^\hush,k)$ is not
a well-defined category ($\Hom$-sets might be large). A convenient
solution is to only consider functors that commute with filtered
colimits. Then each such functor from say $k\amod$ to $k\amod$ is
completely determined by its restriction to the full subcategory
spanned by finitely generated projective $k$-modules, and since this
category is small, the problem does not arise. The same works for
complexes (and the subcategory of perfect complexes). In our
example, all large monoidal categories will be of this sort, so we
will adopt this point of view. By abuse of notation, we will still
use notation of the form $\Fun(\C^\hush,k)$ for the category for
functors that commute with filtered colimits.
\end{remark}

Let us now construct the cyclic power trace functor \eqref{F.cycl}
more canonically. Fix a monoidal category $\C$ and an integer $l
\geq 1$, and define a category $\C^\hush_l$ by the cartesian square
\begin{equation}\label{pi.l.c}
\begin{CD}
\C^\hush_l @>{\pi_l}>> \C^\hush\\
@V{\rho_l}VV @VV{\rho}V\\
\Lambda_l @>{\pi_l}>> \Lambda.
\end{CD}
\end{equation}
Then the functor $i_l:\Lambda_l \to \Lambda$ fits into a commutative
square
\begin{equation}\label{i.l.c}
\begin{CD}
\C^\hush_l @>{i_l}>> \C^\hush\\
@V{\rho_l}VV @VV{\rho}V\\
\Lambda_l @>{i_l}>> \Lambda.
\end{CD}
\end{equation}
Explicitly, $\C^\hush_l$ is the category of pairs $\langle
[n],\{c_v\} \rangle$, $[n] \in \Lambda_l$, $c_v \in \C$, $v \in
V([n])$. Then the top arrow in \eqref{i.l.c} is the functor
sending a sequence $\{c_v\}$, $v \in V([n])$ to the same sequence
repeated $l$ times. In particular, it sends cartesian maps to
cartesian maps.

We now observe that for any trace functor $F$ from $\C$ to the
category of $k$-modules, with the corresponding object $F^\hush \in
\Fun(\C^\hush.k)$, the objects $\pi_{l!}i_l^*F^\hush$,
$\pi_{l*}i_l^*F^\hush$ also correspond to trace functors from $\C$
to $k\amod$. If $\C=k\amod$ and $F$ is the tautological functor,
$\pi_{l!}i_l^*F^\hush$ corresponds to the cyclic power trace functor
\eqref{F.cycl}.

\subsection{Quotients of the conjugate filtration.}

We now fix a perfect field $k$ of characteristic $p$, and we let
$\C=C_\idot(k)$ be the category of complexes of $k$-vector
spaces. Moreover, denote by $I_\idot$ the complex in
$\Fun(\C^\hush,k)$ corresponding to the identity trace functor
$C_\idot(k) \to C_\idot(k)$, and denote by $I^{[p]}_\idot$ the
complex $I_\idot$ equipped with the $p$-the rescaling of the stupid
filtration.

Then the projection $\pi_p:\C^\hush_p \to \C^\hush$ of
\eqref{pi.l.c} is a bifibration with fiber $\ppt_p$, so that we can
consider relative Tate cohomology functor $\pi_p^\flat$. By base
change \cite[Lemma 1.7]{ka0}, we have $\rho^* \circ \pi_p^\flat
\cong \pi_p^\flat \circ \rho^*$, so that \eqref{sum.pi} yields a
direct sum decomposition
\begin{equation}\label{sum.pi.hush}
\pi_p^\flat k \cong \rho^*\K_\idot(k)[1]((u)) = \bigoplus_{i \in
  \Z}\rho^*\K_\idot(k)[2i+1]
\end{equation}
in the derived category $\D(\C^\hush,k)$. We also have the
multiplication map \eqref{mult.pi} and all that it entails --- in
particular, the isomorphisms \eqref{i.n.mult} for complexes
$E_\idot$ in $\Fun(\C^\hush_p,k)$ that are tight in the obvious
sense. We note that the pullback $i_p^*I_\idot$ of the tautological
complex $I_\idot$ is tight.

Now fix some resolution data for $\pi_p$, so that $\pi_p^\flat$ is
defined as a complex, and consider the complex
$$
C_\idot = \tau^F_{[0,1]}\pi_p^\flat i_p^*I^{[p]}_\idot
$$
in the category $\Fun(\C^\hush,k)$. It must correspond to some trace
functor from $\C=C_\idot(k)$ to itself. Explicitly, the trace
functor sends a complex $V_\idot$ to the complex
\begin{equation}\label{c.v.d}
C_\idot(V_\idot) = \tau^F_{[0,1]}\vC^\hdot(\Z/p\Z,V^{\otimes_k
  p}_\idot),
\end{equation}
where we equip $V^{\otimes_k p}_\idot$ with the $p$-th rescaling of
the stupid filtration. By Proposition~\ref{ten.prop}, we have a
natural sequence
\begin{equation}\label{a.b}
\begin{CD}
0 @>>> V^{(1)}_\idot[1] @>{b}>> C_\idot(V_\idot) @>{a}>> V^{(1)}_\idot
@>>> 0
\end{CD}
\end{equation}
of functorial complexes of $k$-vector spaces that is {\em
  quasiexact} in the sense of \cite[Definition 1.2]{ka5} --- this
means that $a \circ b = 0$, the map $a$ is surjective, the map $b$
is injective, and the complex $\Ker a/\Im b$ is acyclic. The map $a$
corresponds to a map
$$
a:C_\idot \to I_\idot^{(1)}
$$
in the category $\Fun(\C^\hush,k)$, where $I^{(1)}_\idot$ is the
Frobenius twist of the tautological complex $I_\idot$.

Moreover, consider the ring $W_2(k)$ of second Witt vectors of the
field $k$, and let $\C_1=C_\idot(W_2(k))$ be the category of
complexes of flat $W_2(k)$-modules. Denote by $q:\C_1 \to \C$ the
quotient functor sending a complex $V_\idot$ to its quotient
$V_\idot/p$. Note that $\C_1$ is a monoidal category, and $q$ is a
monoidal functor. Moreover, extend \eqref{c.v.d} to $W_2(k)$-modules
by setting
\begin{equation}\label{c.v.d.1}
C_\idot(V_\idot) =
\tau^F_{[0,1]}\vC^\hdot(\Z/p\Z,V^{\otimes_{W_2(k)} p}_\idot)
\end{equation}
for any complex $V_\idot \in \C_1$, where as in \eqref{c.v.d}, we
equip $V^{\otimes_{W_2(k)} p}_\idot$ with the $p$-th rescaling of
the stupid filtration. Then we have the following somewhat
surprising result.

\begin{lemma}\label{lft.le}
\begin{enumerate}
\item For any complex $V_\idot \in \C_1$, we have a short exact
  sequence of complexes
$$
\begin{CD}
0 @>>> C_\idot(V_\idot/p) @>{p}>> C_\idot(V_\idot) @>{q}>>
C_\idot(V_\idot/p) @>>> 0,
\end{CD}
$$
where $p$ stands for multiplication by $p$, and $q$ is the quotient
map.
\item For any $V_\idot \in \C$, let $\overline{C}_\idot(V_\idot)$ be
  the kernel of the map $a$ of \eqref{a.b}, and for any $V_\idot \in
  \C_1$, let $\wt{C}_\idot(V_\idot) =
  C_\idot(V_\idot)/p\overline{C}_\idot(V_\idot/p)$. Then the
  composition map
$$
\begin{CD}
\wt{C}_\idot(V_\idot) @>{q}>> C_\idot(V_\idot/p) @>{a}>> V_\idot^{(1)}
\end{CD}
$$
is a quasiisomorphism, and the functor $\wt{C}_\idot:\C_1 \to \C_1$
factors through the quotient functor $q:\C_1 \to \C$.
\end{enumerate}
\end{lemma}

\proof{} \thetag{i} is \cite[Lemma 6.9]{ka5}, and \thetag{ii} is
\cite[Proposition 6.11]{ka5}.
\endproof

We note that $\wt{C}_\idot(-)$ is the cokernel of a map of trace
functors, thus itself inherits the structure of a trace
functor. Then Lemma~\ref{lft.le}~\thetag{ii} implies that this trace
functor is actually defined on $\C$ --- namely, we have the
following.

\begin{corr}\label{W.corr}
Let $\wt{C}_\idot$ be the complex in $\Fun(\C_1^\hush,W_2(k))$
corresponding to the trace functor $\wt{C}_\idot(-)$ of
Lemma~\ref{lft.le}~\thetag{ii}. Then there exists a complex $W_\idot
\in \Fun(\C,W_2(k))$ such that
$$
\wt{C}_\idot \cong q^*W_\idot,
$$
where $q:\C_1^\hush \to \C^\hush$ is induced by the monoidal
quotient functor $q:\C_1 \to \C$.
\end{corr}

\proof{} The quotient functor $q$ is surjective on isomorphism
classes of objects, so that it suffices to prove that the action of
morphisms in $\C_1^\hush$ on $\wt{C}_\idot$ canonically factors
through $q$. For cartesian morphisms, this is clear, so it suffices
to check it for morphisms in the categories $\C_1^V([n])$, $[n] \in
\Lambda$. This immediately follows from
Lemma~\ref{lft.le}~\thetag{ii}.
\endproof

\proof[Proof of Proposition~\ref{char.2.prop}.] Every $k$-vector
space can be considered as a $W_2(k)$-module via the quotient map
$W_2(k)$, so that we have a natural functor from $k$-vector sapces
to $W_2(k)$-module. Moreover, for any small category $I$, we can
apply this pointwise and obtain a functor
$$
\xi:\Fun(I,k) \to \Fun(I,W_2(k)).
$$
This functor is exact and fully faitful. It induces a functor
$\xi:\D(I,k) \to \D(I,W_2(k))$, and for any complex $E_\idot$ in
$\Fun(I,k)$, we have a natural quasiisomorphism
$$
H_\idot(I,\xi(I_\idot)) \cong \xi(H_\idot(I,E_\idot)).
$$
In particular, we can take $I=\Lambda$. Then by \eqref{hc.hla} and
\eqref{hc.hh}, to construct an isomorphism \eqref{char.2}, it
suffices to construct a functorial isomorphism
$$
\xi(\gr^0_V\pi_{p\flat}i_p^*A^\hush_\idot) \cong
\xi(\K(A^{(1)\hush}_\idot))
$$
in the derived category $\D(\Lambda,W_2(k))$.

Consider first the universal situation. Denote
$$
R_\idot = \tau^F_{[1,2]}\pi_p^\flat i_p^*I^{[p]}_\idot[-1],
$$
and recall that we have the multiplication map \eqref{mult.pi} for
the bifibration $\pi_p:\C^\hush_p \to \C^\hush$. This map is
compatible with filtered truncations, so by \eqref{sum.pi.hush}, it
yields a map
$$
\rho^*\K_\idot(k) \otimes_k C_\idot \cong \tau_{[1,2]}\pi_p^\flat
k[-1] \otimes_k \tau^F_{[0,1]}\pi_p^\flat i_p^*I^{[p]}_\idot \to
\tau_{[1,2]}\pi_p^\flat i_p^*I^{[p]}_\idot[-1] = R_\idot.
$$
This map induces a map
$$
\rho^*\K_\idot \otimes \xi(C_\idot) \to
\xi(R_\idot).
$$
On the other hand, we have the complex $W_\idot$ of
Corollary~\ref{W.corr} and the natural map $q:W_\idot \to
\xi(C_\idot)$. Denote by $b$ the composition map
$$
\begin{CD}
\rho^*\K_\idot \otimes W_\idot @>{\id \otimes q}>> \rho^*\K_\idot
\otimes \xi(C_\idot) @>>> \xi(R_\idot).
\end{CD}
$$
Then $b$ is filtered, hence compatible with the filtration
$\tau^\hdot$. Moreover, since the composition $W_\idot \to
\xi(C_\idot) \to \xi(I^{(1)}_\idot)$ is a quasiisomorphism by
Lemma~\ref{lft.le}, the source and the target of the map $b$ only
have two non-trivial associated graded quotients $\gr^i_\tau$, for
$i=0$ and $1$, and in both cases, $\gr^i_\tau(b)$ is one of the
isomorphisms \eqref{i.n.mult} for the tight complex
$I^{[p]}_\idot$. Therefore $b$ is an isomorphism in the derived
category $\D(\C^\hush,W_2(k))$.

Now let $\alpha:\Lambda \to \C^\hush$ be the section of the
projection $\rho:\C^\hush \to \Lambda$ corresponding to the DG
algebra $A_\idot$, and consider the induced isomorphism
$$
\alpha^*(b):\K_\idot(\alpha^*W_\idot) \cong \alpha^*(\rho^*\K_\idot \otimes
W_\idot) \to \alpha^*\xi(R_\idot).
$$
Then the right-hand side is exactly
$\xi(\gr^0_V\pi_{p\flat}i_p^*A^\hush_\idot)$, and the left-hand side is
naturally isomorphic to $\K_\idot(\alpha^*\xi(I^{(1)}_\idot)) \cong
\K_\idot(\xi(A^{(1)\hush}_\idot)) \cong
\xi(\K_\idot(A^{(1)\hush}_\idot))$.
\endproof

\section{Degeneration.}\label{dege.sec}

We now turn to degeneration results for the spectral sequences for
cyclic homology. There are two statements: one for the conjugate
spectral sequence \eqref{conj.sp}, and one for the Hodge-to-de Rham
spectral sequence \eqref{hdr.sp}.

\subsection{Conjugate degeneration.}

Recall that for any field $k$, a {\em square-zero extension}
$A'_\idot$ of a DG algebra $A_\idot$ over $k$ by an
$A_\idot$-bimodule $M_\idot$ is a filtered DG algebra $\langle
A'_\idot,\tau^\hdot \rangle$ such that $\tau^0A'_\idot = A'_\idot$,
$\tau^2A'_\idot = 0$, $\gr^0_\tau A'_\idot \cong A_\idot$, and $\gr^1_\tau
A'_\idot$ is quasiisomorphic to $M_\idot$ as a bimodule over
$\gr^0_\tau A'_\idot \cong A_\idot$. Recall also that up to a
quasiisomorphisms, square-zero extensions are classified by elements
in the {\em reduced Hochschild cohomology group}
$$
\overline{HH}^2(A_\idot,M_\idot) = \Ext^1_{A_\idot^o \otimes
  A_\idot}(I_\idot,M_\idot),
$$
where $I_\idot$ is the kernel of the augmentation map $A^o_\idot
\otimes_k A_\idot \to A_\idot$ (see e.g. \cite[Subsection 4.3]{ka0}
but the claim is completely standard). Reduced Hochschild cohomology
groups are related to the usual ones by the long exact sequence
\begin{equation}\label{red.tri}
\begin{CD}
\overline{HH}^\hdot(A_\idot,M_\idot) @>>> HH^\hdot(A_\idot,M_\idot)
@>>> M_\idot @>>>
\end{CD}
\end{equation}
In particular, if $\overline{HH}^2(A_\idot,M_\idot)=0$ for some
$M_\idot$, then every square-zero extension $A'_\idot$ of $A_\idot$
by $M_\idot$ splits --- there exists a DG algebra $A''_\idot$ and a
map $A''_\idot \to A'_\idot$ such that the composition map
$A''_\idot \to A'_\idot \to A_\idot$ is a quasiisomorphism.

\medskip

Now fix a perfect field $k$ of some positive characteristic $p =
\cchar k$, and assume given a DG algebra $A_\idot$ over $k$.

\begin{theorem}\label{conj.thm}
Assume that the DG algebra $A_\idot$ over the field $k$ satisfies
the following two properties:
\begin{enumerate}
\item There exist a DG algebra $\wt{A}_\idot$ over the second Witt
  vectors ring $W_2(k)$ and a quasiisomorphism $\wt{A}_\idot
  \lotimes_{W_2(k)} k \cong A_\idot$.
\item The reduced Hochschild cohomology $\overline{HH}^i(A_\idot)$
  vanishes for $i \geq 2p$.
\end{enumerate}
Then the conjugate spectral sequence \eqref{conj.sp} degenerates at
first trem, so that there exists an isomorphism $\bHP_\idot(A_\idot)
\cong HH_\idot(A^{(1)}_\idot)((u^{-1}))$.
\end{theorem}

\begin{remark}
By Proposition~\ref{char.2.prop}, the spectral sequence
\eqref{conj.sp} also exists for $p=2$, and Theorem~\ref{conj.thm}
holds in this case, too. Of course in this case, the condition
\thetag{ii} is pretty strong.
\end{remark}

\proof{} By Proposition~\ref{P.dege.prop}, it suffices to construct
a filtered DG algebra $A'_\idot$ over $k$ and a filtered map
$A'_\idot \to P_\idot(A_\idot)$ such that the composition map
$A'_\idot \to A^{(1)}_\idot$ is a filtered quasiisomorphism (the
assumption ``$p$ is odd'' in Proposition~\ref{P.dege.prop} is only
needed to insure that the conjugate spectral sequence is
well-defined, and it is not used in the proof).

To define the DG algebra $P_\idot(A_\idot)$, we need to choose
multiplicative resolution data for the group $G=\Z/p\Z$, and we are
free to do it in any way we like. Note that the $G$-action on the DG
algebra $A^{\otimes_k p}_\idot$ extends to the action of the
symmetric group $\Sigma_p$. In particular, we have an action of the
normalizer $\wh{G} = (\Z/p\Z) \rtimes (\Z/p\Z)^* \subset \Sigma_p$
of $G \subset \Sigma_p$. Choose some multiplicative resolution data
for $\wh{G}$, and restrict it to $G \subset G$. Then the resulting
DG algebra $P_\idot(A_\idot)$ carries a natural action of
$(\Z/p\Z)^* = \wh{G}/G$. This action preserves the filtration
$\tau^\hdot$, and the augmentation map
$$
a:P_\idot(A_\idot) \to \gr^0_\tau P_\idot(A_\idot) \cong A^{(1)}_\idot
$$
is $(\Z/p\Z)^*$-invariant. For $i \geq 1$, we have the isomorphisms
\eqref{gr.P} induced by the isomorphisms \eqref{i.n}; invariantly,
they can be written as
$$
\gr^i_\tau P_\idot(A_\idot) \cong A^{(1)}_\idot[i] \otimes_k
\vH^{-i}(G,k),
$$
where the group $(\Z/p\Z)^*$ acts on the right-hand side via its
action on the Tate comology group $\vH^{-i}(G,k)$.

Now, the cohomology $H^\hdot(G,k)$ is given by \eqref{c.coho}. The
group $(\Z/p\Z)^*$ acts trivially on the generator $\eps$, and it
acts on the generator $u$ via its standard one-dimensional
representation given by the action on $\Z/p\Z \subset k$. Therefore we
have
$$
\vH^{-i}(G,k)^{(\Z/p\Z)^*} \cong k
$$
if $i=0,1 \mod 2(p-1)$, and $0$ otherwise. We conclude that
if we denote
$$
\overline{P}_\idot(A_\idot) = P_\idot(A_\idot)^{(\Z/p\Z)^*} \subset
P_\idot(A_\idot)
$$
and equip this DG algebra with the filtration induced by
$\tau^\hdot$, then we have
\begin{equation}\label{gr.P.b}
\gr^i_\tau \overline{P}_\idot(A_\idot) \cong
\begin{cases}
A^{(1)}_\idot[i], &\quad i=0,1 \mod 2(p-1),\\
0, &\quad \text{otherwise}.
\end{cases}
\end{equation}
If we denote by $e:\overline{P}_\idot(A_\idot) \to P_\idot(A_\idot)$
the embedding map, then $\gr^i_\tau(e)$ is a quasiisomorphism for
$i=0,1 \mod 2p$, and $\gr^i_\tau(e)=0$ otherwise.

To prove degeneration, it suffices to construct a filtered DG
algebra $A'_\idot$ and a filtered map $s:A'_\idot \to
\overline{P}_\idot(A_\idot)$, since then we can simply compose it
with the embedding $e$. Moreover, a filtered algebra $\langle
A'_\idot,\tau^\hdot \rangle$ is completely defined by its quotients
$A^n_\idot = A'_\idot/\tau^{n+1}A'_\idot$, $n \geq 0$, together with
the quotient maps $r_n:A^{n+1}_\idot \to A^n_\idot$, and by
assumption, $A^0_\idot$ must be identified with the Frobenius twist
$A^{(1)}_\idot$ of the DG algebra $A_\idot$. Thus if we denote
$$
P^n_\idot =
\overline{P}_\idot(A_\idot)/\tau^{n+1}\overline{P}_\idot(A_\idot),
\qquad n \geq 1,
$$
and denote by $p_n:P^{n+1}_\idot \to P^n_\idot$ the quotient maps,
then it suffices to construct a collection of DG algebras
$A^n_\idot$ over $k$ for all $n \geq 1$, equipped with DG algebra
maps $r_n:A^{n+1}_\idot \to A^n_\idot$, $s_n:A^n_\idot \to
P^n_\idot$ such that
\begin{itemize}
\item for every $n \geq 1$, we have $s_n \circ r_n = p_n \circ
  s_{n+1}$, and the composition map $a \circ e \circ s_n:A^n_\idot
  \to A^{(1)}_\idot$ is a quasiisomorphism.
\end{itemize}
We use induction on $n$. To start it, we take $n=1$; we need to find
a DG algebra $A^1_\idot$ over $k$ and a map $s_1:A^1_\idot \to
P^1_\idot$ such that $a \circ e \circ s_1:A^1_\idot \to
A^{(1)}_\idot$ is a quasiisomorphism. By \cite[Proposition
  6.13]{ka5}, this is possible precisely because the DG algebra
$A_\idot$ satisfies the assumption \thetag{i} of the Theorem.

For the induction step, assume given $A^{n-1}_\idot$ and
$s_{n-1}:A^{n-1}_\idot \to P^{n-1}_\idot$, and consider the DG
algebra $A''_\idot$ defined by the cartesian square
$$
\begin{CD}
A''_\idot @>>> A^{n-1}_\idot\\
@VVV @VV{s_{n-1}}V\\
P^n_\idot @>{p_{n-1}}>> P^{n-1}_\idot.
\end{CD}
$$
Then up to a quasiisomorphism, $A''_\idot$ is a square-zero
extension of $A^{(1)}_\idot$ by
$\gr^n_\tau\overline{P}_\idot(A_\idot)$, and finding $A^n_\idot$
with the maps $s_n$, $r_{n-1}$ satisfying \thetag{$\bullet$} is
equivalent to finding a DG algebra $A^n_\idot$ and a map $A^n_\idot
\to A''_\idot$ such that the composition map
$$
\begin{CD}
A^n_\idot @>>> A''_\idot @>>> A^{(1)}_\idot
\end{CD}
$$
is a quasiisomorphism. In other words, we have to split the
extension $A''_\idot$. The obstruction to doing this lies in the
reduced Hochschild cohomology group
$\overline{HH}^2(A^{(1)}_\idot,\gr^n_\tau\overline{P}_\idot(A_\idot))$. By
\eqref{gr.P.b}, this group vanishes unless $i=0,1 \mod 2p$, and in
this case, we have
$$
\overline{HH}^2(A^{(1)}_\idot,\gr^n_\tau\overline{P}_\idot(A_\idot))
\cong \overline{HH}^{2+n}(A^{(1)}_\idot).
$$
Since $n \geq 2$ and $n=0,1 \mod 2(p-1)$, we have $n \geq 2(p-1)$,
and then this reduced Hochschild cohomology group vanishes by the
assumption \thetag{ii}.
\endproof

\begin{remark}
The condition \thetag{ii} of Theorem~\ref{conj.thm} is slightly
unnatural: while Hochschild homology and the conjugate spectral
sequence are derived Mori\-ta-invariant, reduced Hochschild cohomology
groups are not (because of the third term in \eqref{red.tri}). One
would like to have the same statement but with $HH^\hdot(-)$ instead
of $\overline{HH}^\hdot(-)$. The simplest way to obtain such a
statement would be to repeat the whole argument for DG categories
instead of DG algebras. In fact, \cite{ka4} also deals with the DG
category case, so that this looks like a straightforward
exercise. However, since our main interest is in degeneration in
$\cchar 0$, we do not go into it to save space.
\end{remark}

\subsection{Hodge-to de Rham degeneration.}

We now fix a field $K$ of characteristic $0$, and a DG algebra
$A_\idot$ over $K$. In this case, a Hodge-to-de Rham degeneration
theorem is an immediate corollary of Theorem~\ref{conj.sp}, and the
argument is exactly the same as in \cite[Subsection 5.3]{ka0}. We
reproduce it for the sake of completeness and for the convenience of
the reader.

\begin{theorem}\label{hdr.thm}
Assume that the DG algebra $A_\idot$ is homologically smooth and
homologically proper. Then the Hodge-to-de Rham spectral sequence
\eqref{hdr.sp} degenerates, so that there exists an isomorphism
$HP_\idot(A_\idot) \cong HH_\idot(A_\idot)((u))$.
\end{theorem}

We recall that {\em homologically proper} simply means that
$A_\idot$ is a perfect complex over $K$ (in particular, it is
homologically bounded).

\proof{} By a theorem of B.To\"en \cite{to}, there exists a finitely
generated subring $R \subset K$ and a homologically smooth and
homologically proper DG algebra $A^R_\idot$ over $R$ such that
$A_\idot \cong A^R_\idot \lotimes_R K$. Since $R$ is finitely
generated, the residue field $k=R/\m$ for any maximal ideal $\m
\subset R$ is a finite, hence perfect field of some characteristic
$p$.  Since $A^R_\idot$ is homologically proper and homologically
smooth, there is at most a finite number of non-trivial Hochschild
homology groups $HH_\idot(A^R_\idot)$ and reduced Hochschild
cohomology groups $\overline{HH}^\hdot(A^R_\idot)$, and these groups
are finitely generated $R$-modules. Then there exists a constant $N$
such that $\overline{HH}^i(A^R_\idot)=0$ for $i \geq N$. Moreover,
localizing $R$ if necessary, we can further assume that
$HH_i(A^R_\idot)$ is a projective finitely generated $R$-module for
every $i$, and that for any maximal ideal $\m \subset R$, $p =
\cchar R/\m$ is non-trivial in $\m/\m^2$ (that is, $p$ is unramified
in $R$), and $2p > N$. Then for any $\m \subset R$ with $k = R/\m$,
the DG algebra $A^k_\idot = A^R_\idot \lotimes_R k$ satisfies the
assumptions of Theorem~\ref{conj.sp}. Therefore we have an
isomorphism
$$
\bHP_\idot(A^k_\idot) \cong HH_\idot(A^k_\idot)^{(1)}((u^{-1}))
$$
of finite-dimensional graded $k$-vector spaces. Since
$HH_\idot(A^k_\idot)$ is concentrated in a finite range of degrees,
we can replace Laurent power series in $u^{-1}$ with Laurent power
series in $u$, and since $A^R_\idot$ is homologically smooth and
proper, $A^k_\idot$ is also homologically smooth and proper. In
particular, it is cohomologically bounded, so that
Theorem~\ref{per.thm}~\thetag{ii} allows us to replace
$\bHP_\idot(A^k_\idot)$ with $HP_\idot(A^k_\idot)$. We thus have an
isomorphism
$$
HP_\idot(A^k_\idot) \cong HH_\idot(A^k_\idot)^{(1)}((u))
$$
of finite-dimensional graded $k$-vector spaces, so that the
Hodge-to-de Rham spectral sequence \eqref{hdr.sp} for $A^k_\idot$
degenerates for dimension reasons.

Finally, since all the Hochschild homology $R$-modules
$HH_i(A^R_\idot)$ are finitely generated and projective, and any
differential in the Hodge-to-de Rham spectral sequence for
$A^R_\idot$ vanishes modulo any maximal ideal $\m \subset R$, the
differential must vanish identically. Thus the Hodge-to-de Rham
spectral sequence for the DG algebra $A^R_\idot$ degenerates, and
then so does the Hodge-to-de Rham spectral sequence for $A_\idot =
A^R_\idot \lotimes_R K$.
\endproof

{\footnotesize

\bigskip

\noindent
{\sc
Steklov Math Institute, Algebraic Geometry section\\
\mbox{}\hspace{30mm}and\\
Center for Geometry and Physics, IBS, Pohang, Rep. of Korea
}

\medskip

\noindent
{\em E-mail address\/}: {\tt kaledin@mi.ras.ru}
}

\end{document}